\documentclass[a4paper,12pt,reqno]{amsart}
\usepackage{lmodern,amsmath,amssymb,mathtools,amsthm}
\usepackage{enumerate}
\usepackage{graphics,graphicx}
\usepackage{stackengine}
\usepackage{bigints}
\usepackage{geometry}
\usepackage{multirow}
\usepackage{wasysym}
\usepackage{hyperref}
\usepackage{floatrow,morefloats}
\usepackage{epsfig}
\usepackage{epstopdf}
\usepackage{longtable}
\usepackage{caption,subcaption}
\newtheorem{definition}{Definition}
\newtheorem{lemma}{Lemma}
\newtheorem{theorem}{Theorem}
\newtheorem{corollary}{Corollary}
\newtheorem{remark}{Remark}

\newtheorem{problem}{Problem}

\newtheorem{proposition}{Proposition}
\numberwithin{equation}{section}
\numberwithin{definition}{section}
\numberwithin{lemma}{section}
\numberwithin{theorem}{section}
\numberwithin{corollary}{section}
\numberwithin{proposition}{section}
\numberwithin{example}{section}
\numberwithin{remark}{section}

\usepackage{blindtext}

\title{Orthogonality of a new family of $q$-Sobolev type polynomials}

\author{Neha$^\dagger$ and A. Swaminathan $^{\#\ddagger}$}
\address{$^\dagger$Department of Mathematics, Indian institute of techonology, Roorkee-247667, Uttarakhand-India}
\email{neharani0777@gmail.com, neha@ma.iitr.ac.in}
\address{$^\ddagger$Department of Mathematics, Indian institute of techonology, Roorkee-247667, Uttarakhand-India}
\email{mathswami@gmail.com, a.swaminathan@ma.iitr.ac.in}
\allowdisplaybreaks
\bigskip
\newmuskip\pFqmuskip

\newcommand*\pFqskip{8mu}
\catcode`,\active
\newcommand*\pFq{\begingroup
	\catcode`\,\active
	\def ,{\mskip\pFqskip\relax}%
	\dopFq
}
\catcode`\,12
\def\dopFq#1#2#3#4#5{%
	{}_{#1}F_{#2}\biggl[\genfrac..{0pt}{}{#3}{#4};#5\biggr]%
	\endgroup
}

\begin{document}

	%%%%%%%%%%%%%%%%%%%% ABSTRACT %%%%%%%%%%%%%%%%%%%%%%%%
	\begin{abstract}
		In this work, we introduce and construct specific $q$-polynomials that are desired from the well-established families of $q$-orthogonal polynomials, namely little $q$-Jacobi polynomials and $q$-Laguerre polynomials, respectively. We examine these newly constructed $q$-polynomials and observe that they possess integral representations of little $q$-Jacobi polynomials and $q$-Laguerre polynomials. These polynomials solve a third-order $q$-difference equation and display an unconventional four-term recurrence relation. This unique recurrence relation makes us categorize them as $q$-Sobolev-type orthogonal polynomials. This motivation leads to defining the general Sobolev-type orthogonality for $q$-polynomials. Special cases of these polynomials are also explored and discussed. Furthermore, we delve into the behavior of these $q$-orthogonal polynomials of Sobolev type as the parameters approach $1$. We also examine their zeros and interlacing properties.
	\end{abstract}
	%%%%%%%%%%%%%%%%%%%%%%%%%%%%%%%%%%%%%%%%%%%%%%%%%%%%%%%
	
	\keywords{q-Orthogonal polynomials; q-Laguerre polynomials; little q-Jacobi polynomials; basic hypergeometric function; q-Sobolev orthogonal polynomials; zeros.\\
		\# Corresponding author}
	
	\maketitle
	
	\markboth{Neha and A. Swaminathan}{Special functions}

	%%%%%%%%%%%%%%%%%%%%%%%%%%%%%%%%%% SECTION INTRODUCTION %%%%%%%%%%%%%%%%%%%%%%%%%%%%%%%%%%%%%%%%%%%%%
	
	\section{Introduction}
	\label{Introduction}
	In the realm of $q$-special functions (\cite{E. Koelink_q-special functions,1960_Rainville_special functions}), $q$-polynomials assume a fundamental and essential role, making significant contributions to a diverse array of problems. These problems encompass a wide range of fields, including Eulerian series, continued fractions \cite{S. Morier-Genoud_$q$-deformed rationals and $q$-continued fractions}, $q$-algebra, quantum groups, and $q$-oscillators, among others. Within the context of $q$-special functions, our primary focus is placed on the $q$-orthogonal polynomials \cite{1990_Gasper_Basic hypergeometric series} and these polynomials solve specific second-order $q$-difference equations  of the form \cite{2010_Koekoek Lesky Swarttow_Hypergeometric orthogonal polynomials},
	\begin{align}
		&(wx^2+2bx+v)(\mathcal{D}_q^2p_n)(x)+(2\epsilon x+\eta)(\mathcal{D}_qp_n)(x)=\frac{[n]_q}{q^n}(w[n-1]
		_q+2\epsilon)p_n(qx),\\
		&n\in \mathbb{N}_0=\mathbb{N}\cup\{0\},\quad w,b,v,\epsilon,\eta\in\mathbb{R},\quad 0<q<1, \label{q-difference equation of order two} \notag
	\end{align}
	where the $q$-difference operator and $q$-number, respectively, are given as
	\begin{align*}
		\mathcal{D}_qh(x):=\left\{\begin{array}{ll}
			\dfrac{h(qx)-h(x)}{(q-1)x}, & \mbox{$x\neq0$}\\
			h'(0), & \mbox{$x=0.$}
		\end{array}\right. \quad\mbox{and}\quad
		[a]_q=\frac{q^a-1}{q-1}.
	\end{align*}
	We define 
	\begin{align*}
		\mathcal{D}_q^0h:=h\quad \mbox{and}\quad\mathcal{D}_q^ih:=\mathcal{D}_q(\mathcal{D}_q^{i-1}h),\quad i=1,2,....
	\end{align*}
	If the function $h(x)$ is differentiable at $x$, it is straightforward to demonstrate that
	\begin{align*}
		\lim_{q \to 1}\mathcal{D}_qh(x)=h'(x).
	\end{align*}
	
	The significance of $q$-orthogonal polynomials lies in their wide-ranging applications across various mathematical disciplines. They play a crucial role in generalizing classical orthogonal polynomials (\cite{1959_Szego_orthogonal polynomials}), and are deeply connected to algebraic structures, special functions, and various areas of theoretical physics (see \cite{2012_Andrews_Ramanujan identities,2003_Kemp_Characterizations_discrete distributions,1992_Koelink_Koornwinder_q-special functions,2019_Sri Ranga_Complementary {R}omanovski-{R}outh polynomials}). These polynomials are expressed through a basic hypergeometric series, commonly known as the Heine series, with the detailed information provided in \cite{1869_Thomae_Heine Series}. Depending on the parameter $q$, they are alternatively referred to as basic (or $q$-) hypergeometric series. Heine used the notation ${_m\phi_p}$ to represent these series (see \cite{1990_Gasper_Basic hypergeometric series}).
	
	Moreover, within the domain of $q$-orthogonal polynomials \cite{1990_Gasper_Basic hypergeometric series}, our specific focus centers on $q$-Sobolev orthogonal polynomials. It is important to highlight that Sobolev orthogonal polynomials have undergone extensive exploration in recent years. Althammer made the initial strides in understanding Sobolev orthogonal polynomials, drawing inspiration from Lewis' work \cite{1947_Lewis_Polynomial least square approximations}. In \cite{1962_Althammer}, Althammer demonstrated that the polynomials, orthogonal under the inner product defined below, constitute Sobolev orthogonal polynomials.
	\begin{equation}
		\begin{aligned}
			\langle l,\,m\rangle_{\lambda}=\int_{-1}^{1}l(z)m(z)dz+\lambda\int_{-1}^{1}l'(z)m'(z)dz,\quad \lambda>0.\label{definition of Sobolev orthogonal polynomials}
		\end{aligned}
	\end{equation}
	
	In \cite{1962_Althammer}, Legendre-Sobolev polynomials, a subset of Sobolev orthogonal polynomials, were examined. Orthogonal polynomials with inner products involving difference operators were explored in a series of papers \cite{1996_Bavinck_charlier polynomials,1997_Bavinck_linear perturbation}. These polynomials, orthogonal with respect to inner products involving difference operators, are termed Sobolev orthogonal polynomials.
	
	In \cite{2020_Sergey_Classical type SOP}, the discussion centers around classical-type Sobolev orthogonal polynomials. However, our present research aims to expand upon these polynomials, resulting in the identification of $q$-Sobolev orthogonal polynomials of classical type. These $q$-Sobolev orthogonal polynomials serve as a more comprehensive framework encompassing classical-type Sobolev orthogonal polynomials, which include well-known instances such as Jacobi, Laguerre, and Chebyshev polynomials (\cite{2018_Duran_Differential equations for discrete Jacobi-Sobolev orthogonal polynomials,1991_koekoek_generalized laguerre polynomials,2000_Koekoek_jacobi polynomials,2017_Manas_Asymptotics for varying discrete Sobolev orthogonal polynomials,2021_Markett_Sobolev orthogonal polynomials of Laguerre- and Jacobi-type}). By doing so, they provide a broader framework that includes classical-type Sobolev orthogonal polynomials as special cases. This generalization allows for a deeper understanding of the relationships between families of orthogonal polynomials. Recent developments in Sobolev orthogonal polynomials are presented in \cite{2019_Sri Ranga_{S}obolev orthogonal polynomials on the unit circle,2023_Ranga_{S}obolev-type orthogonal polynomials
		on the unit circle based on a {$q$}-difference operator,2023_Sri Ranga_Coherent pairs of measures of the second kind on the real line}.
	
	In this work, we maintain the condition of $0 < q < 1$ consistently. Let $\{a_i\}_{i=1}^{m}$ and $\{b_j\}_{j=1}^{p}$ be the complex numbers such that $b_j$ is not equal to $q^{-n}$ with $n\in \mathbb{N}_0$,\quad $j=1,2,\ldots p.$
	The $q$-hypergeometric series, $_m\phi_p$ with variable $z$ is defined as follows
	\[{_m\phi_p}\left[{_{b_1,\ldots,b_p}^{a_1,\ldots,a_m};q;z}\right] := \sum\limits_{k = 0}^\infty  {\frac{{{{(a_1;q)}_k}\ldots{{(a_m;q)}_k}}}{{{{(q;q)}_k}{{(b_1;q)}_k}\ldots{{(b_p;q)}_k}}}}\big[(-1)^kq^{\frac{k(k-1)}{2}}\big]^{(1+p-m)}z^k,\]
	where
	\begin{align*}
		{{(b;q)}_k}=\left\{\begin{array}{cl}
			1, & \mbox{$k=0$,} \\
			(1-b)(1-bq)\ldots(1-bq^{k-1}), & \mbox{$k=1,2,\ldots,$}
		\end{array}\right.
	\end{align*}
	is the {\it{q-shifted factorial}}.\\
	The $q$-hypergeometric series can be understood as a $q$-analogue of the hypergeometric series, and this relationship can be described as follows
	\begin{align*}
		\lim_{q \to 1}{_m\phi_p}\left[{_{q^{b_1},\ldots,q^{b_p}}^{q^{a_1},\ldots,q^{b_m}};q;(q-1)^{1+p-m}z}\right]&=\pFq{m}{p}{a_1,\ldots,a_m}{b_1,\ldots,b_p}{z}.
	\end{align*}
	The hypergeometric function is expressed as
	\begin{align*}
		\pFq{m}{p}{a_1,\ldots,a_m}{b_1,\ldots,b_p}{z}
		:=\sum_{k=0}^{\infty} \frac{(a_1)_k\ldots(a_m)_k}{(b_1)_k\ldots(b_p)_k}\frac{z^k}{k!},
	\end{align*}
	and $(s)_l$ is the rising factorial given by
	\begin{align*}
		(s)_l=\left\{\begin{array}{cl}
			1,  & \mbox{$l=0$,} \\
			s(s+1)\ldots(s+l-1), & \mbox{$l=1,2,\ldots.$}
		\end{array}\right.
	\end{align*}
	{\bf Motivation for the problem}: We specifically focus on two types of $q$-orthogonal polynomials: namely, $q$-Laguerre and little $q$-Jacobi polynomials.
	\begin{itemize}
		\item The 	Little $q$-Jacobi polynomials are given by 
		\begin{align}\label{little q-jacobi}
			\mathcal{P}_n^{(\gamma,\xi)}(z;q)={_2\phi_1}\left[{_{q^{\gamma+1}}^{q^{-n},q^{n+\gamma+\xi+1}};q;qz}\right],\quad \gamma,\xi>-1, \;n \in \mathbb{N}_0.
		\end{align}
		\item The $q$-Laguerre polynomials are given by
		\begin{align}\label{q-laguerre}
			\mathcal{L}_n^{(\gamma)}(z;q)={_1\phi_1}\left[{_{q^{\gamma+1}}^{q^{-n}};q;-q^{n+\gamma+1}z}\right],\quad \gamma>-1, \;n \in \mathbb{N}_0.
		\end{align}
	\end{itemize}
	The specific case of little $q$-Jacobi polynomials corresponds to little $q$-Legendre polynomials when $\gamma=\xi=0$. We can extend these polynomials by including additional parameters in the denominator and  numerator of the basic hypergeometric functions, defining them as
	\begin{itemize}
		\item Generalized little $q$-Jacobi polynomials are defined as
		\begin{align}
			\label{Generalized little q-Jacobi polynomials}
			\mathfrak{P}_n^{(\gamma,\xi)}(z,c;q):={_3\phi_2}\left[{_{q^{\gamma+1},q^{c+1}}^{q^{-n},q^{n+\gamma+\xi+1},q};q;qz}\right],\quad \gamma,\xi>-1, \;n \in \mathbb{N}_0.
		\end{align}
		\item Generalized $q$-Laguerre polynomials are defined as
		\begin{align}\label{Generalized q-laguerre polynomials}
			\mathfrak{L}_n^{(\gamma)}(z,c;q):={_2\phi_2}\left[{_{q^{\gamma+1},q^{c+1}}^{q^{-n},q};q;-q^{n+\gamma+1}z}\right],\quad \gamma>-1, \;n \in \mathbb{N}_0.
		\end{align}
	\end{itemize}
	When $c=0$, these generalized polynomials reduce to the original little $q$-Jacobi and $q$-Laguerre polynomials, as follows\\
	$\mathfrak{P}_n^{(\gamma,\xi)}(z,0;q)=\mathcal{P}_n^{(\gamma,\xi)}(z;q)$ and $\mathfrak{L}_n^{(\gamma)}(z,0;q)=\mathcal{L}_n^{(\gamma)}(z;q).$
	In our analysis, we consider $c$ to be a positive number, providing insight into the nature of these generalized polynomials. Understanding the orthogonality properties of the original $q$-Laguerre and little $q$-Jacobi polynomials, coupled with knowledge about the structure of their generalized counterparts, sets the groundwork. By establishing a relationship between these generalized and original polynomials and utilizing the known orthogonality properties of the original ones, we can deduce the orthogonality properties of the generalized little $q$-Jacobi and generalized $q$-Laguerre polynomials. This connection is essential in unraveling the properties of the generalized polynomials based on the known properties of the original ones. To bridge this connection, we can use $q$-difference operator, a powerful tool that allows us to reduce generalized $q$-orthogonal polynomials to their original $q$-orthogonal polynomial forms.
	
	The structure of this work is as follows: In Section~\ref{Hypergeometric $q$-Sobolev orthogonal polynomials}, we present integral representations of $\mathfrak{L}^{(\gamma)}_n(z,c;q)$ and $\mathfrak{P}_n^{(\gamma,\xi)}(z,c;q)$ through propositions~\ref{the integral representation of generalized little q-Jacobi polynomials} and \ref{the integral representation of generalized q-laguerre polynomials}, establishing their connection to $q$-Laguerre and Little $q$-Jacobi polynomials, respectively. We also establish that these polynomials adhere to a $q$-difference equation of a hypergeometric type (refer to propositions~\ref{relation between OP and Sobolev type OP for little q-Jacobi}, \ref{relation between OP and Sobolev type OP for q-laguerre}). Furthermore, we provide linear recurrence relations for these polynomials, as demonstrated in Theorems~\ref{theorem for recurrence relation for generalized little q-jacobi polynomials} and \ref{theorem for recurrence relation for generalized q-laguerre polynomials}. Subsequently, we proceed to characterize the polynomials $\mathfrak{L}^{(\gamma)}_n(z,c;q)$ and $\mathfrak{P}_n^{(\gamma,\xi)}(z,c;q)$, firmly establishing their classification as $q$-Sobolev type orthogonal polynomials. We also discuss all these properties for the special cases of these polynomials. Lastly, in Section \ref{behavior of zeros}, we investigate the behavior of zeros as parameters $\gamma$, $\xi$, and $c$ approach $1$. Additionally, we explore the interlacing properties of zeros for $\mathfrak{L}_n^{(\gamma)}(z,c;q)$ within the range $-1<\gamma<1$ while $c=1$.
	\section{Characterization of generalized $q$-orthogonal polynomials}{\label{Hypergeometric $q$-Sobolev orthogonal polynomials}}
	Consider a non-trivial positive measure on the real line given by $d_q\mu(z)=w_q(z)d_qz$. The system of polynomials $\{p_n(z;q)\}_{n=0}^{\infty}$ with $\text{degree}(p_n(z;q))=n$ is termed $q$-orthogonal on the compact subset $K$ of $\mathbb{R}$ if  
	\begin{align}\label{definition of orthogonality-continuous}
		\langle p_n(z;q),\,p_m(z;q)\rangle=\int_{K}p_n(z;q)p_m(z;q)d_q\mu(z)=A_{n,q}\delta _{mn}, \; A_{n,q}>0, \, n,m \in \mathbb{N}_0,
	\end{align} 
	holds.
	
	The Jacobi polynomials are expressed as
	\begin{align}\label{classical jacobi polynomials}
		P_n(z,\gamma,\xi)=\pFq{2}{1}{-n,n+\gamma+\xi+1}{\gamma+1}{z},\; \gamma,\xi>-1,\;n\in\mathbb{N}_0.
	\end{align}
	The Jacobi polynomials, defined in \eqref{classical jacobi polynomials}, possess a $q$-analogue, which is given by \eqref{little q-jacobi}.
	By introducing the parameter $c$, we extend the Jacobi polynomials, resulting in
	\begin{align*}
		\mathcal{P}_n(z,\gamma,\xi,c)=\pFq{3}{2}{-n,n+\gamma+\xi+1,1}{\gamma+1,c+1}{z},\; \gamma,\xi>-1,c>0,\;n\in\mathbb{N}_0.
	\end{align*}
	The $q$-analogue of these modified polynomials is given by \eqref{Generalized little q-Jacobi polynomials}.
	
	The Laguerre polynomials are expressed as
	\begin{align}\label{laguerre polynomials}
		L_n(z,\gamma)=\pFq{1}{1}{-n}{\gamma+1}{z},\; \gamma>-1,\;n\in\mathbb{N}_0.
	\end{align}
	The Laguerre polynomials have a $q$-analogue, expressed as \eqref{q-laguerre}. The Laguerre polynomials, extended by introducing the parameter $c$ in \eqref{laguerre polynomials}, are given by
	\begin{align*}
		\mathcal{L}_n(z,\gamma,c)=\pFq{2}{2}{-n,1}{\gamma+1,c+1}{z},\;\gamma>-1,\;c>0,\;n\in\mathbb{N}_0.
	\end{align*}
	The generalized Laguerre polynomials have a $q$-analogue, given by \eqref{Generalized q-laguerre polynomials}.
	
	The properties of the polynomials $\mathfrak{P}_n^{(\gamma,\xi)}(z,c;q)$ and $\mathfrak{L}_n^{(\gamma)}(z,c;q)$ (see \cite{2018_Costas-Santos_Analytic properties of some basic HSTOP,1998_Koekoek_laguerre polynomials,1981_Moak_q-analogue of laguerre polynomials}) encompass being solutions to hypergeometric $q$-difference equations and having integral representations  that incorporate little $q$-Jacobi and $q$-Laguerre polynomials. The subsequent propositions summarize these properties. The interpretation of the polynomials $\mathfrak{P}^{(\gamma,\xi)}_n(z,c;q)$ and $\mathfrak{L}^{(\gamma)}_n(z,c;q)$ as $q$-integral representations (see \cite{1910_Jackson_On q-integrals,2005_Sole Kac_On integral representations of q-gamma}) relies on the properties of the $q$-beta, $q$-gamma and basic hypergeometric functions.
	
	The $q$-gamma function $\Gamma_q(e)$, an extension of Euler's gamma function in the $q$-analogue framework, was introduced in \cite{1869_Thomae_Heine Series}. Subsequently, it was further examined in \cite{1910_Jackson_On q-integrals} as an infinite product, and its expression is given by
	\begin{align*}
		\Gamma_q(e)=\frac{(1-q)_q^{e-1}}{(1-q)^{e-1}},\quad e>0,
	\end{align*}
	where
	\begin{align*}
		(x+y)^n_q=\prod_{j=0}^{n-1}(x+q^jy),\quad n\in\mathbb{N}.
	\end{align*}
	The $q$-beta function $\beta_q(a,b)$ is given by the formula
	\begin{align*}
		\beta_q(a,b)=\frac{\Gamma_q(a)\Gamma_q(b)}{\Gamma_q(a+b)},
	\end{align*}
	and it also possesses the $q$–integral representation, provided below.
	\begin{align*}
		\beta_q(a,b)=\int_{0}^{1}t^{a-1}(1-qt)_q^{b-1}d_qt,\quad Re(a),\,Re(b)>0.
	\end{align*}
	The expression of the $q$-beta function will play a vital role in the context of propositions~\ref{the integral representation of generalized little q-Jacobi polynomials} and \ref{the integral representation of generalized q-laguerre polynomials}.
	\begin{proposition}{\label{the integral representation of generalized little q-Jacobi polynomials}}
		Polynomials $\mathfrak{P}_n^{(\gamma,\xi)}(z,c;q)$, with $\gamma,\xi>-1$ and $c>0$, are represented by the following integral expression
		\begin{align*}
			\mathfrak{P}_n^{(\gamma,\xi)}(z,c;q)=\frac{1-q^{c}}{1-q}\int_{0}^{1}(1-qt)_q^{c-1}\mathcal{P}_n^{(\gamma,\xi)}(zt;q)d_qt,\quad |z|<1.
		\end{align*}
		\begin{proof}
			We know that
			\begin{align*}
				\mathfrak{P}_n^{(\gamma,\xi)}(z,c;q)&=\sum_{i=0}^{n}\frac{(q^{-n};q)_i(q^{n+\gamma+\xi+1};q)_i(q;q)_i}{(q^{\gamma+1};q)_i(q^{c+1};q)_i}\frac{(qz)^i}{(q;q)_i}.
			\end{align*}
			Since 
			\begin{align*}
				\frac{(q;q)_i}{(q^{c+1};q)_i}&=\frac{1}{\binom{i+c}{i}_q}=\frac{\Gamma_q{(i+1)}\Gamma_q{(c+1)}}{\Gamma_q{(i+c+1)}}=\frac{1-q^{c}}{1-q}\frac{\Gamma_q{(i+1)}\Gamma_q{(c)}}{\Gamma_q{(i+c+1)}}\\
				&=\frac{1-q^{c}}{1-q}\beta_q{(i+1,c)}=\frac{1-q^{c}}{1-q}\int_{0}^{1}t^i(1-qt)_q^{c-1}d_qt.
			\end{align*}
			Then 
			\begin{align*}
				\mathfrak{P}_n^{(\gamma,\xi)}(z,c;q)&=\sum_{i=0}^{n}\frac{(q^{-n};q)_i(q^{n+\gamma+\xi+1};q)_i}{(q^{\gamma+1};q)_i}\Bigg(\frac{1-q^{c}}{1-q}\int_{0}^{1}t^i(1-qt)_q^{c-1}d_qt\Bigg)\frac{(qz)^i}{(q;q)_i}\\
				&=\frac{1-q^{c}}{1-q}\int_{0}^{1}\sum_{i=0}^{n}\frac{(q^{-n};q)_i(q^{n+\gamma+\xi+1};q)_i}{(q^{\gamma+1};q)_i}\frac{q^i}{(q;q)_i}(zt)^i(1-qt)_q^{c-1}d_qt\\
				&=\frac{1-q^{c}}{1-q}\int_{0}^{1}(1-qt)_q^{c-1}\mathcal{P}_n^{(\gamma,\xi)}(zt;q)d_qt.
			\end{align*}
			Hence the result.
		\end{proof}
	\end{proposition}
	
	To derive the integral representation for the polynomials $\mathfrak{L}_n^{(\gamma)}(z,c;q)$, we can utilize an useful property, avoiding the need to repeat the process used for deriving the integral representation for $\mathfrak{P}_n^{(\gamma,\xi)}(z,c;q)$.
	\begin{align}\label{limiting value of generalized little q-jacobi polynomials}
		\mathfrak{L}_n^{(\gamma)}(z,c;q)=\lim_{b \to -\infty}\mathfrak{P}_n^{\big(\gamma,\frac{\log b}{\log q}\big)}\bigg(-\frac{z}{bq},c;q\bigg), \;n\in\mathbb{N}_0,\;z\in\mathbb{R},\;\gamma>-1,\;c\geq0.
	\end{align}
	In Proposition~\ref{the integral representation of generalized little q-Jacobi polynomials}, substitute $z$ with $-\frac{z}{bq}$, and then take the limit as $b$ approaches $-\infty$, to establish the following proposition.
	\begin{proposition}{\label{the integral representation of generalized q-laguerre polynomials}}
		Polynomials $\mathfrak{L}_n^{(\gamma)}(z,c;q)$, with $\gamma>-1$ and $c>0$, are represented by the following integral expression
		\begin{align*}
			\mathfrak{L}_n^{(\gamma)}(z,c;q)=\frac{1-q^{c}}{1-q}\int_{0}^{1}(1-qt)_q^{c-1}\mathcal{L}_n^{(\gamma)}(zt;q)d_qt, \quad z\in\mathbb{C}.
		\end{align*}
	\end{proposition}
	In an upcoming result, we establish that the basic hypergeometric function adheres to a specific $q$-difference equation \cite{1990_Gasper_Basic hypergeometric series} of a particular type:
	\begin{align}\label{hypergeometric type q-difference equation}
		\big(\Delta\Delta_{b_1/q}\ldots\Delta_{b_p/q}\big)v(z;q)=z\big(\Delta_{a_1}\ldots\Delta_{a_m}\big)v(zq^{1+p-m};q),
	\end{align}
	where
	\begin{align*}
		v(z;q)={_m\phi_p}(a_1,\ldots,a_m;b_1,\ldots,b_p;q;z)={_m\phi_p}\left[{_{b_1,\ldots,b_p}^{a_1,\ldots,a_m};q;z}\right],
	\end{align*}
	and 
	\begin{equation*}
		\begin{aligned}
			\Delta_{b}f(x)=bf(qx)-f(x).
		\end{aligned}
	\end{equation*}
	The equation \eqref{hypergeometric type q-difference equation} can be seen as the $q$-analogue of the equation:
	\begin{equation*}
		\begin{aligned}
			\bigg(\theta\prod_{j=1}^{s}(\theta+b_j-1)-z\prod_{i=1}^{r}(\theta+a_i)\bigg)w=0,
		\end{aligned}
	\end{equation*}
	where $\theta=z\frac{d}{dz}$ and it is satisfied by ${}_r F_s$ (a generalized hypergeometric function) (\cite{1960_Rainville_special functions}, p. no. 73).
	
	\begin{proposition}{\label{theorem of higher order q-difference equation for little q-Jacobi polynomials}}
		The polynomials $y(z;q)=\mathfrak{P}_n^{(\gamma,\xi)}(z,c;q),\,(\gamma,\xi>-1,\;c>0)$ satisfy the following q-difference equation:
		\begin{equation*}
			\begin{aligned}
				&(q^{\gamma+c}-q^{\gamma+\xi+2}z)q^3(q-1)^2z^2(\mathcal{D}_q^3y)(z;q)+\big(q^{\gamma+c}(q+1)-q^{\gamma}-q^{c}+q^{\gamma+\xi+1}z+q^{1-n}z\\
				&+q^{n+\gamma+\xi+2}z-q^{\gamma+\xi+2}(q^2+q+1)z\big)q(q-1)z(\mathcal{D}_q^2y)(z;q)+\big((q^{\gamma+1}-1)(q^{c+1}-1)\\
				&-(q^{n+\gamma+\xi+2}+q^{n+\gamma+\xi+1}-q^{n+\gamma+\xi}-1)q(q^{-n+1}-1)z-q^{-n}(q-1)(q^{n+\gamma+\xi+1}-1)z\big)\\
				&(\mathcal{D}_qy)(z;q)-(q^{-n}-1)(q^{n+\gamma+\xi+1}-1)y(z;q)=0.
			\end{aligned}
		\end{equation*}
		\begin{proof}
			Since the basic hypergeometric function satisfies a $q$-difference equation \eqref{hypergeometric type q-difference equation}, the polynomials $\mathfrak{P}_n^{(\gamma,\xi)}(z,c;q)$, with $(\gamma,\xi>-1,\;c>0)$, also satisfy the following $q$-difference equation
			\begin{align*}
				(\Delta\Delta_{q^{\gamma}}\Delta_{q^{c}})y(z;q)
				=z(\Delta_{q^{-n}}\Delta_{q^{n+\gamma+\xi+1}}\Delta_q)y(z;q),
			\end{align*}
			where $\Delta_ay(z;q)=a(q-1)z\mathcal{D}_qy(z;q)+(a-1)y(z;q)$. After simplifications, we get the desired result.
		\end{proof}
	\end{proposition}
	\begin{remark}
		It is not feasible to derive the $q$-difference equation for $\mathfrak{L}_n^{(\gamma)}(z,c;q)$ by simply applying the limit as demonstrated in \eqref{limiting value of generalized little q-jacobi polynomials}.
	\end{remark}
	\begin{proposition}{\label{theorem of higher order q-difference equation for q-laguerre polynomials}}
		The polynomials $y(z;q)=\mathfrak{L}_n^{(\gamma)}(z,c;q),\,(\gamma>-1,\,c>0)$ satisfy the following q-difference equation:
		\begin{equation*}
			\begin{aligned}
				&(q^{\gamma+c}-zq^{-n+1})q^3(q-1)^2z^2(\mathcal{D}^3_qy)(z;q)+\big(q^{\gamma+c}(q+1)-q^{c}-q^{\gamma}-q^{-n+1}(q^2+q+1)z\\
				&+(q^{-n}+q)z\big)q(q-1)z(\mathcal{D}^2_qy)(z;q)+\big((q^{\gamma+1}-1)(q^{c+1}-1)-(q^2+q-1)(q^{-n+1}-1)z\\
				&-q^{-n}(q-1)z\big)(\mathcal{D}_qy)(z;q)-(q^{-n}-1)y(z;q)=0.
			\end{aligned}
		\end{equation*}
		\begin{proof}
			We can establish this similarly to Proposition~\ref{theorem of higher order q-difference equation for little q-Jacobi polynomials} by computing
			\begin{align*}
				(\Delta\Delta_{q^{\gamma}}\Delta_{q^{c}})y(z;q)
				=z(\Delta_{q^{-n}}\Delta_q)y(qz;q).
			\end{align*}
			Hence, we have obtained the desired result.
		\end{proof}
	\end{proposition}
	The $q$-difference equation presented in Proposition~\ref{theorem of higher order q-difference equation for little q-Jacobi polynomials}  resembles the generic form of a $q$-difference equation, expressed as
	\begin{align}
		L_{\alpha}y(z;q)+\lambda\mathfrak{D}_q^{t}y(z;q)=0,\label{strumn-liouville type q-equation}
	\end{align} 
	where
	\begin{equation*}
		\begin{aligned}
			L_{\alpha}y(z;q)&=\sum_{i=0}^{\alpha}l_i(z;q)\mathcal{D}_q^iy(z;q).
		\end{aligned}
	\end{equation*}
	\begin{remark}
		Propositions~\ref{theorem of higher order q-difference equation for little q-Jacobi polynomials} and \ref{theorem of higher order q-difference equation for q-laguerre polynomials} demonstrate that $\mathfrak{P}_n^{(\gamma,\xi)}(z,c;q)$ and $\mathfrak{L}_n^{(\gamma)}(z,c;q)$ conform to a q-difference equation given by \eqref{strumn-liouville type q-equation} with specific parameters: $\alpha=3$ and $t=0$ for $\mathfrak{P}_n^{(\gamma,\xi)}(z,c;q)$, and $\alpha=3$ and $t=0$ for $\mathfrak{L}_n^{(\gamma)}(z,c;q)$.
	\end{remark}
	Now, let us move forward to establish a recurrence relation \cite{2023_Paco_Higher-order recurrence relations} for the polynomials  $\mathfrak{P}^{(\gamma,\xi)}_n(z,c;q)$ and $\mathfrak{L}^{(\gamma)}_n(z,c;q)$. We state this as the following theorems. 
	
	\begin{theorem}{\label{theorem for recurrence relation for generalized little q-jacobi polynomials}}
		The polynomials $y_n(z;q)=\mathfrak{P}_n^{(\gamma,\xi)}(z,c;q) \,(\gamma,\xi>-1,c>0)$ exhibit the following four term recurrence relation:
		\begin{equation}\label{recurrence relation for genralized little q-jacobi polynomials}
			\begin{aligned}
				&\mu_{q,1}(n)y_{n-2}(z;q)+(\mu_{q,2}(n)+x\mu_{q,5}(n))y_{n-1}(z;q)+(\mu_{q,3}(n)\\
				&+x\mu_{q,6}(n))y_n(z;q)+\mu_{q,4}(n)y_{n+1}(z;q)=0,\quad n \in \mathbb{N}_0,
			\end{aligned}
		\end{equation}
		where $\mu_{q,j}(n)\;(1\leq j\leq6)$ for  $n\geq2$ is given by the relations \eqref{Coefficients of recurrence realtions}. The coefficients $\mu_{q,1}(n),\mu_{q,4}(n),\mu_{q,5}(n) \;and\;\mu_{q,6}(n)$ for $n=0,1$ are defined by the same expression.
		Moreover,
		%\vspace{-0.7em}
		\begin{equation*}
			\begin{aligned}
				\mu_{q,2}(0)&=0,\\
				\mu_{q,3}(0)&=(-1)\big((1-q^{\gamma+1})(1-q^{c+1})[q^{\gamma+\xi+1};q]_2+(1-q^{\gamma})(1-q^{c})[q^{\gamma+\xi+2};q]_2\\
				&-(1-q^{\gamma})(1-q^{c})[q^{\gamma+\xi+2};q]_2\big),\\
				\mu_{q,2}(1)&=\frac{q(1-q^{\gamma+1})(1-q^{c+1})(1-q^{\gamma+\xi})[q^{\gamma+\xi+4};q]_2}{(1-q^{\gamma+\xi+2})}+q(1-q^{\gamma})(1-q^{c})[q^{\gamma+\xi+4};q]_2\\
				&-q(1-q^{\gamma})(1-q^{c})[q^{\gamma+\xi+4};q]_2-q(1-q^{\gamma+2})(1-q^{c+2})[q^{\gamma+\xi+2};q]_2\\
				&+(1+q)(1-q^{\gamma+2})(1-q^{c+2})(1-q^{\gamma+\xi+1})(1-q^{\gamma+\xi+3})\\
				&-\frac{(1+q)(1-q^{\gamma+1})(1-q^{c+1})(1-q^{\gamma+\xi+1})[q^{\gamma+\xi+4};q]_2}{(1-q^{\gamma+\xi+2})},\\
				\mu_{q,3}(1)&=(-1)[q^{\gamma+\xi+1};q][q^{\gamma+\xi+3};q]\Big(\frac{q(1-q^{\gamma+1})(1-q^{c+1})(1-q^{\gamma+\xi})(1-q^{\gamma+\xi+4})}{[q^{\gamma+\xi+2};q]_2}\\
				&+(1+q)(1-q^{\gamma+2})(1-q^{c+2})-\frac{(1+q)(1-q^{\gamma+1})(1-q^{c+1})(1-q^{\gamma+\xi+4})}{(1-q^{\gamma+\xi+2})}\Big),
			\end{aligned}
		\end{equation*}
		$\mathfrak{P}_{-1}^{(\gamma,\xi)}(z,c;q)=0$ \quad\mbox{and} \quad $\mathfrak{P}_{-2}^{(\gamma,\xi)}(z,c;q)=0.$
		%	\vspace{-1em}
		\begin{proof}
			The statement can be directly verified for $n=0$ and $n=1$. Proof of the case $n\geq2$ differs from the discussions in the sequel. It is given separately in Section~\ref{proof of theorem} at the end of the manuscript.
		\end{proof}
	\end{theorem}
	To obtain the recurrence relation for the polynomials $\mathfrak{L}_n^{(\gamma)}(z,c;q)$, we employ \eqref{limiting value of generalized little q-jacobi polynomials} in the recurrence relation for $\mathfrak{P}_n^{(\gamma,\xi)}(z,c;q)$.
	%	\vspace{-1em}
	\begin{theorem}{\label{theorem for recurrence relation for generalized q-laguerre polynomials}}
		Polynomials $y_n(z;q)=\mathfrak{L}_n^{(\gamma)}(z,c;q)$ exhibit the following recurrence relation:
		\begin{equation}
			\begin{aligned}
				&q\mu'_{q,1}(n)y_{n-2}(z;q)+(q\mu'_{q,2}(n)+x\mu'_{q,5}(n))y_{n-1}(z;q)+(q\mu'_{q,3}(n)\\
				&+x\mu'_{q,6}(n))y_n(z;q)+q\mu'_{q,4}(n)y_{n+1}(z;q)=0,\quad n \in \mathbb{N}_0,\label{recurrence relation for generalized q-laguerre polynomials}
			\end{aligned}
		\end{equation}
		where 
		\vspace{-0.5em}
		\begin{equation*}
			\begin{aligned}
				\mu'_{q,3}(n)&=(-1)q^{3n+2\gamma+1}\Bigg(\frac{(1-q^{n+1})(1-q^{\gamma+n+1})(1-q^{c+n+1})}{(1-q)}+q^2(1-q^{\gamma+n})(1-q^{c+n})\\
				&-\frac{q^2(1-q^{n+1})(1-q^{\gamma+n})(1-q^{c+n})}{(1-q)}\Bigg),\\
				\mu'_{q,4}(n)&=q^{3n+2\gamma+1}(1-q^{\gamma+n+1})(1-q^{c+n+1}),\\
				\mu'_{q,5}(n)&=q^{5n+3\gamma+3}(1-q^n),\\
				\mu'_{q,6}(n)&=-q^{5n+3\gamma+3}(1-q^{n+1}), \\
				\mu'_{q,2}(n)&=\frac{q^{3n+2\gamma+1}(1-q^{\gamma+n+1})(1-q^{c+n+1})[q^{n+1};q]_2}{(1-q)^2}\\
				&+\frac{q^{3n+2\gamma+3}(1-q^{\gamma+n})(1-q^{c+n})(1-q^n)}{(1-q)}\\
				&-\frac{q^{3n+2\gamma+3}(1-q^{\gamma+n})(1-q^{c+n})[q^{n+1};q]_2}{(1-q)^2}\\
				&-\frac{q^{3n+2\gamma+1}(1-q^{\gamma+n+1})(1-q^{c+n+1})[q^{n+1};q]_2}{[q^2;q]_2}\\
				&-\frac{q^{3n+2\gamma+5}(1-q^{\gamma+n-1})(1-q^{c+n-1})(1-q^n)}{(1-q)}\\
				&+\frac{q^{3n+2\gamma+5}(1-q^{\gamma+n-1})(1-q^{c+n-1})[q^{n+1};q]_2}{[q^2;q]_2},\\
				{\mbox{and}}\quad	\mu'_{q,1}(n)&=-\mu'_{q,2}(n)-\mu'_{q,3}(n)-\mu'_{q,4}(n),\;n\geq2.
			\end{aligned}
		\end{equation*}
		The coefficients $\mu'_{q,3}(n),\mu'_{q,4}(n),\;and\;\mu'_{q,6}(n)$ for $n=0,1$ are defined by the same expressions.
		Moreover,
		\begin{equation*}
			\begin{aligned}
				&\mu'_{q,5}(0)=0,\;\mu'_{q,2}(0)=0,\;\mu'_{q,1}(0)=0,\;\mu'_{q,1}
				(1)=0,\;\mu'_{q,5}(1)=q^{3\gamma+8}(1-q),\\
				&\mu'_{q,2}(1)=q^{2\gamma+5}(1-q^{\gamma+2})(1-q^{c+2})-q^{2\gamma+7}(1-q^{\gamma+1})(1-q^{c+1}),\\
				&\mathfrak{L}_{-1}^{(\gamma)}(z,c;q)=\mathfrak{L}_{-2}^{(\gamma)}(z,c;q)=0.
			\end{aligned}
		\end{equation*}
		\vspace{-1.2em}
		\begin{proof}
			Write the expression \eqref{recurrence relation for genralized little q-jacobi polynomials} for $\mathfrak{P}_n^{\big(\gamma,\frac{logb}{logq}\big)}\Big(-\frac{z}{bq},c;q\Big)$. Then dividing it by $b^3$ and using $b \to-\infty$ gives the desired result.
		\end{proof}
	\end{theorem}
	\begin{remark}
		If a system of polynomials follow a TTRR, Favard's theorem \cite{1978_Chihara_orthogonal polynomials} confirms their orthogonality. However, when dealing with the polynomials $\mathfrak{P}_n^{(\gamma,\xi)}(z,c;q)$ and $\mathfrak{L}_n^{(\gamma)}(z,c;q)$, which adhere to a four-term recurrence relation, we cannot definitively assert their orthogonality.
	\end{remark}
	To evaluate whether these polynomials are orthogonal, consider the $q$-difference operator defined below
	\begin{equation}\label{q-difference equation of Sobolev OP-1}
		\begin{aligned}
			\mathfrak{D}_q^{t}y(z;q)&=\sum_{j=0}^{t}d_j(z;q)\mathcal{D}_q^jy(z;q),\quad t\in\mathbb{N}_0,
		\end{aligned}
	\end{equation}
	where $d_j(z;q)$ are polynomials of variable $z$. By employing this operator, we can simplify $\mathfrak{P}_n^{(\gamma,\xi)}(z,c;q)$ polynomials and $\mathfrak{L}_n^{(\gamma)}(z,c;q)$ polynomials to their classical versions, specifically the $\mathcal{P}_n^{(\gamma,\xi)}(z;q)$ polynomials and $\mathcal{L}_n^{(\gamma)}(z;q)$ polynomials. This simplification is denoted by the expression
	\begin{align}
		\mathfrak{D}_q^{t} y_n(z;q)=p_n(z;q), \label{q-difference equation with required condition}
	\end{align}
	where $p_n(z;q)$ represents the little $q$-Jacobi polynomials and $q$-Laguerre polynomials that are $q$-orthogonal. The accuracy and suitability of the equation \eqref{q-difference equation with required condition} have been proven and confirmed in Theorems~\ref{relation between OP and Sobolev type OP for little q-Jacobi} and \ref{relation between OP and Sobolev type OP for q-laguerre}. However, it is important to note that this simplification does not apply universally to all generalized $q$-orthogonal polynomials. Hence, our attention is directed exclusively towards $\mathcal{P}_n^{(\gamma,\xi)}(z;q)$ polynomials and $\mathcal{L}_n^{(\gamma)}(z;q)$ polynomials, making them the primary focus. This leads us to the orthogonality relationship implied by equation \eqref{definition of orthogonality-continuous}, as
	\begin{align}	\label{definition of sobolev orthogonality}
		\int_{K}\mathfrak{D}_q^{t} y_n(z;q)\mathfrak{D}_q^{t} y_m(z;q)d_q\mu(x)=A_{n,q}\delta _{mn}, \; A_{n,q}>0,\, m,n \in \mathbb{N}_0.
	\end{align}
	
	Given that equation \eqref{definition of sobolev orthogonality} involves $y_n(z;q)$, we can represent it as $\langle y_n, y_m \rangle_S$ for simplicity as
	\begin{align*}
		\langle y_n, y_m \rangle_S=\int_{K}\mathfrak{D}_q^{t} y_n(z;q)\mathfrak{D}_q^{t} y_m(z;q)d_q\mu(x)=A_{n,q}\delta _{mn}.
	\end{align*}
	In the following discussion, we show that the set of polynomials $\{y_n(z;q)\}_{n=0}^{\infty}$ fulfills the $q$-difference equation.

	The polynomials $p_n(z;q)$ satisfy the $q$-difference equation
	\begin{align}
		\sum_{i=0}^{\gamma'}\tilde{l}_i(z;q)\mathcal{D}_q^ip_n(z;q)=\lambda_{n,q}p_n(qz;q), \quad n\in \mathbb{N}_0. \label{q-difference equation of arbitrary higher order}
	\end{align}
	We can denote the operator $\sum_{i=0}^{\gamma'}\tilde{l}_i(z;q)\mathcal{D}_q^i$ as $\tilde{L}_{\gamma'}$. Then, as a consequence of \eqref{q-difference equation with required condition}, the system of polynomials $\{y_n(z;q)\}_{n=0}^{\infty}$ adheres to the subsequent $q$-difference equation
	\begin{align}
		\tilde{L}_{\gamma'}\mathfrak{D}_q^{t}y_n(z;q)=\lambda_{n,q}\mathfrak{D}_q^{t}y_n(z;q), \quad n\in \mathbb{N}_0. \label{q-difference equation of Sobolev OP}
	\end{align}	
	Before proceeding to Theorem~\ref{relation between OP and Sobolev type OP for little q-Jacobi} concerning the well-definedness of \eqref{q-difference equation with required condition}, it is crucial to note that \eqref{definition of orthogonality-continuous} specifically apply to continuous $q$-orthogonal polynomials. Now, we introduce the concept of orthogonality for discrete $q$-orthogonal polynomials (see \cite{2010_Koekoek Lesky Swarttow_Hypergeometric orthogonal polynomials}). If the weight function has a countably infinite number of points of increase, then the weight function (measure) can be written as follows
	\begin{align*}
		d_q\mu(z)=w_q(z)d_qz=\sum_{k=0}^{\infty}\frac{1}{2^k}\delta(z-q^k).
	\end{align*}
	Then \eqref{definition of orthogonality-continuous} implies,
	\begin{align*}
		\langle p_n(z;q),\,p_m(z;q)\rangle=\sum_{k=0}^{\infty}\frac{1}{2^k} p_n(q^k;q)p_m(q^k;q).%=A_{n,q}\delta_{mn}.
	\end{align*}
	Using \eqref{q-difference equation with required condition}, we obtain
	\begin{align}\label{definition of sobolev orthogonality-discrete}
		&\langle y_n(z;q),\,y_m(z;q)\rangle_S=\sum_{k=0}^{\infty}\frac{1}{2^k}\mathfrak{D}_q^{t} y_n(q^k;q)\mathfrak{D}_q^{t} y_m(q^k;q).%=A_{n,q}\delta _{mn}.
	\end{align}
	\begin{theorem}{\label{relation between OP and Sobolev type OP for little q-Jacobi}}
		For arbitrary $c\in\mathbb{N}$, the polynomials $y_n(z;q)=\mathfrak{P}^{(\gamma,\xi)}_n(z,c;q)$ satisfy the orthogonality relation \eqref{definition of sobolev orthogonality-discrete} with $d_q\mu(z)=\sum_{k=0}^{\infty}\frac{(q^{\xi+1};q)_k}{(q;q)_k}(q^{\gamma+1})^k\delta(z-q^k),\; t=c$;\quad$q^{\xi+1}<1$,\;$0<q^{\gamma+1}<1,$
		\begin{align*}
			A_{n,q}&=\Bigg(\frac{(q;q)_{c}}{(1-q)^{c}}\Bigg)^2 \frac{(q^{\gamma+\xi+2};q)_{\infty}}{(q^{\gamma+1};q)_{\infty}}\frac{(1-q^{\gamma+\xi+2})(q^{\gamma+1})^n}{(1-q^{2n+\gamma+\xi+1})}\frac{(q,q^{\xi+1};q)_n}{(q^{\gamma+1},q^{\gamma+\xi+1};q)_n},\quad n\in \mathbb{N}_0,\\
			d_j(z;q)&=\Bigg(\frac{(q;q)_{c}}{(q;q)_j}\Bigg)^2\frac{q^{jc}}{(q;q)_{c-j}(1-q)^{c-j}}z^j, \quad 0\leq j \leq c.
		\end{align*}
		Moreover, the polynomials $\{y_n(z;q)\}_{n=0}^{\infty}$ exhibit the q-difference equation \eqref{q-difference equation of Sobolev OP} with \\$\gamma'=2,\;\tilde{l}_2(z;q)=q^{\gamma+\xi+1}z^2-q^{\gamma-1}
		z,\;\tilde{l}_1(z;q)=\Big(\frac{1-q^{\gamma+\xi+2}}{q(q-1)}z-\frac{1-q^{\gamma+1}}{q^2(q-1)}\Big),\;\tilde{l}_0(z;q)=0\;and\;\\\lambda_{n,q}=\frac{{{[n]_q}}}{{{q^n}}}\Big(q^{\gamma+\xi+1}[n-1]_q+\frac{1-q^{\gamma+\xi+2}}{q(1-q)}\Big).$
		\begin{proof}
			Observe that
			\begin{align*}
				\mathfrak{D}_q^{c}(\mathfrak{P}_n^{(\gamma,\xi)}(z,c;q))&=\mathcal{D}_q^{c}(z^{c}\mathfrak{P}_n^{(\gamma,\xi)}(z,c;q))\\
				&=\mathcal{D}_q^{c}\Bigg(z^{c}\sum_{j=0}^{n}\frac{(q^{-n};q)_j(q^{n+\gamma+\xi+1};q)_j(q;q)_j}{(q^{\gamma+1};q)_j(q^{c+1};q)_j}\frac{(qz)
					^j}{(q;q)_j}\Bigg)\\
				&=\mathcal{D}_q^{c}\Bigg(\sum_{j=0}^{n}\frac{(q^{-n};q)_j(q^{n+\gamma+\xi+1};q)_j}{(q^{\gamma+1};q)_j(q^{c+1};q)_j}q^jz^{j+c}\Bigg)\\
				&=\sum_{j=0}^{n}\frac{(q^{-n};q)_j(q^{n+\gamma+\xi+1};q)_j}{(q^{\gamma+1};q)_j(q^{c+1};q)_j}q^j\mathcal{D}_q^{c}(z^{j+c}).
			\end{align*}
			Clearly
			\begin{align*}
				\mathcal{D}_q^{c}\big(z^{j+c}\big)&=\frac{1}{(1-q)^{c}z^{c}}\sum_{l=0}^{c}(-1)^l\binom{c}{l}_qq^{\binom{l}{2}-(c-1)l}\big(q^lz\big)^{j+c}\\
				&=\frac{1}{(1-q)^{c}}\sum_{l=0}^{c}(-1)^l\binom{c}{l}_qq^{\binom{l}{2}+l+jl}z^{j}\\
				%&=\frac{x^j}{(1-q)^{c}}\sum_{l=0}^{c}(-1)^l\binom{c}{l}_qq^{\binom{l}{2}}q^{(j+1)l}\\
				&=\frac{z^j}{(1-q)^{c}}(q^{j+1};q)_{c}.
			\end{align*}
			Then, we have
			\begin{align*}
				\mathfrak{D}_q^{c}(\mathfrak{P}_n^{(\gamma,\xi)}(z,c;q))&=\sum_{j=0}^{n}\frac{(q^{-n};q)_j(q^{n+\gamma+\xi+1};q)_j}{(q^{\gamma+1};q)_j(q^{c+1};q)_j}q^j\frac{z^j(q^{j+1};q)_{c}}{(1-q)^{c}}\\
				%&=\frac{1}{(1-q)^{c}}\sum_{j=0}^{n}\frac{(q^{-n};q)_j(q^{n+\gamma+\xi+1};q)_j(q;q)_{c}}{(q^{\gamma+1};q)_j(q;q)_j}(qx)^j\\
				&=\frac{(q;q)_{c}}{(1-q)^{c}}\sum_{j=0}^{n}\frac{(q^{-n};q)_j(q^{n+\gamma+\xi+1};q)_j}{(q^{\gamma+1};q)_j}\frac{(qz)^j}{(q;q)_j}
			\end{align*}
			\begin{align}\label{connection between generalized and original polynomials}
				\mathfrak{D}_q^{c}(\mathfrak{P}_n^{(\gamma,\xi)}(z,c;q))&=\frac{(q;q)_{c}}{(1-q)^{c}}\mathcal{P}_n^{(\gamma,\xi)}(z;q).
			\end{align}
			Here, we can choose
			\begin{align*}
				p_n(z;q)=\frac{(q;q)_{c}}{(1-q)^{c}}\mathcal{P}_n^{(\gamma,\xi)}(z;q).
			\end{align*}
			Now, we can express $\mathfrak{D}_q^{t}y(z;q) = \mathcal{D}_q^{c}(z^{c}y(z;q))$ in the form of \eqref{q-difference equation of Sobolev OP-1} as
			\begin{align*}
				\mathfrak{D}_q^{t}y(z;q)&=\mathcal{D}_q^{c}(z^{c}y(z;q))\\
				&=\sum_{j=0}^{c}\binom{c}{j}_q\mathcal{D}_q^{c-j}((q^jz)^{c})\mathcal{D}_q^{j}y(z;q)
				\\
				%&=\sum_{j=0}^{c}\binom{c}{j}_qq^{jc}\frac{(q^{j+1};q)_{c-j}}{(1-q)^{c-j}}x^j\mathcal{D}_q^{j}y(z;q)\\
				&=\sum_{j=0}^{c}\Bigg(\frac{(q;q)_{c}}{(q;q)_j}\Bigg)^2\frac{q^{jc}}{(q;q)_{c-j}(1-q)^{c-j}}z^j\mathcal{D}_q^{j}y(z;q).
			\end{align*}
			Moreover, since $\mathcal{P}_n^{(\gamma,\xi)}(z,c;q)$ is a solution to the $q$-difference equation
			\begin{align*}
				&\frac{{{[n]_q}}}{{{q^n}}}\Bigg(q^{\gamma+\xi+1}[n-1]_q+\frac{1-q^{\gamma+\xi+2}}{q(1-q)}\Bigg)y_n(qz;q)=\\
				&\Bigg(q^{\gamma+\xi+1}z^2-q^{\gamma-1}z\Bigg)(\mathcal{D}_q^2y_n)(z;q)+\Bigg(\frac{1-q^{\gamma+\xi+2}}{q(q-1)}z-\frac{1-q^{\gamma+1}}{q^2(q-1)}\Bigg)
				(\mathcal{D}_qy_n)(z;q),
			\end{align*}
			\eqref{q-difference equation of arbitrary higher order} holds with $\gamma'=2,\,\tilde{l}_2(z;q)=q^{\gamma+\xi+1}z^2-q^{\gamma-1}z,\;\tilde{l}_1(z;q)=\Big(\frac{1-q^{\gamma+\xi+2}}{q(q-1)}z-\frac{1-q^{\gamma+1}}{q^2(q-1)}\Big),\\ \tilde{l}_0(z;q)=0\;and\;\lambda_{n,q}=\frac{{{[n]_q}}}{{{q^n}}}\Big(q^{\gamma+\xi+1}[n-1]_q+\frac{1-q^{\gamma+\xi+2}}{q(1-q)}\Big).$ From the properties of little $q$-Jacobi orthogonal polynomials (see \cite{2010_Koekoek Lesky Swarttow_Hypergeometric orthogonal polynomials}), the orthogonality condition for q-orthogonal polynomials and \eqref{q-difference equation of Sobolev OP}, we obtain the value of $A_{n,q}$ and the $q$-difference equation for $y_n(z;q)$.
		\end{proof}
	\end{theorem}
	The theorem referenced in (\cite{2020_Sergey_Classical type SOP}, Theorem 2) explains the characteristics of the generalization of Jacobi polynomials. We can prove this theorem by examining the limit as $q\to1$ in Theorem~\ref{relation between OP and Sobolev type OP for little q-Jacobi}.
	\begin{corollary}{\rm{\cite{2020_Sergey_Classical type SOP}}}{\label{relation between OP and Sobolev type OP for Jacobi}}
		For arbitrary $c\in\mathbb{N}$, the polynomials
		\begin{align*}
			y_n(z)=\mathcal{P}_n(z,\gamma,\xi,c)=\pFq{3}{2}{-n,n+\gamma+\xi+1,1}{\gamma+1,c+1}{z},\;\gamma,\xi>-1,
		\end{align*}
		satisfy the orthogonality relation \eqref{definition of sobolev orthogonality} with the ordinary differentiation operator and \quad\\ $t=c,\;K=[0,1],\;w(z)=z^{\gamma}(1-z)^{\xi}$,
		\begin{align*}
			A_{n}=\dfrac{\big(\gamma(c+1)\big)^2\gamma{(\gamma+1+n)}\gamma{(\xi+n+1)}}{\bigg(\binom{n+\gamma}{n}\bigg)^2(2n+\gamma+\xi+1)n!\gamma{(n+\gamma+\xi+1)}}, \quad n\in \mathbb{N}_0,
		\end{align*}
		\begin{align*}
			d_j(z)=\Bigg(\frac{c!}{j!}\Bigg)^2\frac{1}{(c-j)!}z^j, \quad 0\leq j \leq c.
		\end{align*}
	\end{corollary}
	The numerical validation of the above corollary, demonstrating that $\mathfrak{P}^{(\gamma,\xi)}_n(z,c;q)$ converges to $\mathcal{P}_n(z,\gamma,\xi,c)$ as $q$ approaches 1, is presented in Table~\ref{zeros of generalized little q-Jacobi polynomials as q approaches to 1}. This table provides evidence of the convergence of zeros for the polynomials $\mathfrak{P}_6^{(0.1,0.2)}(z,1;q)$ as $q\to1$. These zeros match the zeros of the polynomials $\mathcal{P}_6(z,\gamma=0.1,\xi=0.2,c=1)$. The distinct properties regarding the behavior of zeros are elaborated in \cite{2018_Kerstin_Zeros of Jacobi polynomials}, and we delve into these properties in a later section. All the computations were carried out using Mathematica\textsuperscript{\textcopyright}  software on a system equipped with a 12th Gen Intel(R) Core(TM) i5-1235U processor running at 1.30 GHz and 16 GB of RAM.
	
	\begin{table}[ht]
		\begin{center}
			\begin{tabular}{|c c c|c|}
				\hline\hline
				\multicolumn{3}{|c|}{Zeros of $\mathfrak{P}_6^{(0.1,0.2)}(z,1;q)$}&Zeros of\\
				$q=0.99997$&$q=0.99998$&$q=0.999985$&$\mathcal{P}_6(z,\gamma=0.1,\xi=0.2,c=1)$\\
				\hline\hline
				0.0949794&0.0949798&0.090562&0.0949788\\
				\hline
				0.257525&0.257516&0.257931&0.257537\\
				\hline
				0.501793&0.50193&0.496392&0.502\\
				\hline
				0.720398&0.718613&0.720023&0.715285\\
				\hline
				0.886528&0.893004&0.893478&0.90901\\
				\hline
				1.01233&1.00685&1.00458&0.992734\\
				\hline
			\end{tabular}
		\end{center}
		\caption{Zeros of $\mathfrak{P}_6^{(\gamma=0.1,\xi=0.2)}(z,c=1;q)$ when $q \to 1$ converges to zeros of $\mathcal{P}_6(z,\gamma=0.1,\xi=0.2,c=1)$}
		\label{zeros of generalized little q-Jacobi polynomials as q approaches to 1}
	\end{table}
	
	Let us explore the characteristics of generalized $q$-Laguerre polynomials. We relate generalized $q$-Laguerre and $q$-Laguerre polynomials without relying on the limit \eqref{limiting value of generalized little q-jacobi polynomials} in \eqref{connection between generalized and original polynomials} to prevent confusion associated with interchanging the limit with the $q$-difference operator.
	\begin{theorem}{\label{relation between OP and Sobolev type OP for q-laguerre}}
		For arbitrary $c\in\mathbb{N}$, the polynomials $y_n(z;q)=\mathfrak{L}_n^{(\gamma)}(z,c;q)$, $(\gamma>-1)$ satisfy the orthogonality relation \eqref{definition of sobolev orthogonality} with \; $w_q(z)=\dfrac{z^{\gamma}}{(-z;q)_{\infty}},\;t=c,\; K=[0,\infty)$,
		\begin{equation*}
			\begin{aligned}
				A_{n,q}&=\frac{(q^{\gamma+1};q)_n}{(q;q)_nq^n}\frac{(q^{-\gamma};q)_{\infty}}{(q;q)_{\infty}}\Bigg(\frac{(q;q)_{c}}{(1-q)^{c}}\Bigg)^2\Gamma{(-\gamma)}\Gamma{(\gamma+1)
				}, \quad n\in \mathbb{N}_0,\\
				d_j(z;q)&=\Bigg(\frac{(q;q)_{c}}{(q;q)_j}\Bigg)^2\frac{q^{jc}}{(q;q)_{c-j}(1-q)^{c-j}}z^j, \quad 0\leq j \leq c.
			\end{aligned}
		\end{equation*}
		Moreover, the polynomial $y_n(z;q)$ solves the q-difference equation \eqref{q-difference equation of Sobolev OP} with $\gamma'=2,\;\\\tilde{l}_2(z;q)=qz^2+z,\;\tilde{l}_1(z;q)=\frac{q}{q-1}z+\frac{[\gamma+1]_q}{q^{\gamma+1}},\;\tilde{l}_0(z;q)=0\;and\;\lambda_{n,q}=\frac{[n]_q}{q-1}.$
		\begin{proof}
			An easy computation gives
			\begin{align*}
				\mathcal{D}_q^{c}(z^{c}\mathfrak{L}_n^{(\gamma)}(z,c;q))&=\mathcal{D}_q^{c}\Bigg(\sum_{i=0}^{n}\frac{(q^{-n};q)_i(q;q)_i}{(q^{\gamma+1};q)_i(q^{c+1};q)_i}\frac{(-1)^iq^{\binom{i}{2}}(-q^{n+\gamma+1})^i}{(q;q)_i}z^{i+c}\Bigg)\\
				&=\sum_{i=0}^{n}\frac{(q^{-n};q)_i(q;q)_i}{(q^{\gamma+1};q)_i(q^{c+1};q)_i}\frac{(-1)^iq^{\binom{i}{2}}(-q^{n+\gamma+1})^i}{(q;q)_i}\mathcal{D}_q^{c}\big(z^{i+c}\big)\\
				%&=\sum_{j=0}^{n}\frac{(q^{-n};q)_j(q^{j+1};q)_{c}}{(q^{\gamma+1};q)_j(q^{c+1};q)_j}\frac{(-1)^jq^{\binom{j}{2}}(-q^{n+\alpha+1})^j}{(1-q)^{c}}x^{j}\\
				%&=\frac{(q;q)_{c}}{(1-q)^{c}}\sum_{j=0}^{\infty}\frac{(q^{-n};q)_j}{(q^{\alpha+1};q)_j}\frac{(-1)^jq^{\binom{j}{2}}(-q^{n+\alpha+1})^j}{(q;q)_j}x^{j}\\
				&=\frac{(q;q)_{c}}{(1-q)^{c}}\mathcal{L}_n^{(\gamma)}(z;q).
			\end{align*}
			Thus, we can choose
			\begin{align*}
				p_n(z;q)=\frac{(q;q)_{c}}{(1-q)^{c}}\mathcal{L}_n^{(\gamma)}(z;q).
			\end{align*}
			Now, we can express $\mathfrak{D}_q^{t}y(z;q) = \mathcal{D}_q^{c}(z^{c}y(z;q))$ in the form of \eqref{q-difference equation of Sobolev OP-1} as follows
			\begin{align*}
				\mathfrak{D}_q^{t}y(z;q)&=\mathcal{D}_q^{c}(z^{c}y(z;q))\\
				&=\sum_{i=0}^{c}\Bigg(\frac{(q;q)_{c}}{(q;q)_i}\Bigg)^2\frac{q^{ic}}{(q;q)_{c-i}(1-q)^{c-i}}z^i\mathcal{D}_q^{i}y(z;q).
			\end{align*}
			Moreover, since $\mathcal{L}_n^{(\gamma)}(z;q)$ satisfies the $q$-difference equation
			\begin{align*}
				(qz^2+z)(\mathcal{D}_q^2y_n)(z;q)+\Bigg(\dfrac{q}{q-1}z+\dfrac{[\gamma+1]_q}{q^{\gamma+1}}\Bigg)(\mathcal{D}_qy_n)(z;q)=\frac{[n]_q}{q-1}y_n(qz;q),
			\end{align*}
			\eqref{q-difference equation of arbitrary higher order} holds with $\gamma'=2,\,\tilde{l}_2(z;q)=qz^2+z,\,\tilde{l}_1(z;q)=\Big(\frac{q}{q-1}z+\frac{[\gamma+1]_q}{q^{\gamma+1}}\Big),\,\tilde{l}_0(z;q)=0$\;and\\ $\lambda_{n,q}=\frac{[n]_q}{q-1}.$ Upon examining the characteristics of the $q$-Laguerre orthogonal polynomials \cite{2010_Koekoek Lesky Swarttow_Hypergeometric orthogonal polynomials} and considering the relations \eqref{definition of sobolev orthogonality} and \eqref{q-difference equation of Sobolev OP}, we can derive the value of $A_{n,q}$ and the $q$-difference equation for $y_n(z;q).$
		\end{proof}
	\end{theorem}
	The theorem presented in (\cite{2020_Sergey_Classical type SOP}, Theorem 1.) demonstrates the generalization of Laguerre polynomials. We intend to establish the same theorem in the subsequent content through an alternative approach.	
	\begin{corollary}{\rm{\cite{2020_Sergey_Classical type SOP}}}{\label{relation between OP and Sobolev type OP for laguerre}}
		For arbitrary $c\in\mathbb{N}$, the polynomials 
		\begin{align*}
			y_n(z)=\mathcal{L}_n(z,\gamma,c)=\pFq{2}{2}{-n,1}{\gamma+1,c+1}{z},\;\gamma>-1,
		\end{align*}
		satisfy the orthogonality relation \eqref{definition of sobolev orthogonality} with the ordinary differentiation operator and \\ $t=c;\;w(z)=z^{\gamma}e^{-z};\;K=[0,\infty)$, where
		\begin{align*}
			A_n=\dfrac{(\gamma(c+1))^2\gamma{(\gamma+1+n)}}{\bigg(\binom{n+\gamma}{n}\bigg)^2n!}, \quad n\in \mathbb{N}_0
		\end{align*} and
		\begin{align*}
			d_j(z)=\Bigg(\frac{c!}{j!}\Bigg)^2\frac{1}{(c-j)!}z^j, \quad 0\leq j \leq c.
		\end{align*}
		\begin{proof}
			We achieve the desired result by changing $z$ to $(1-q)z$ and then using the limit $q \to 1$ in Theorem~\ref{relation between OP and Sobolev type OP for q-laguerre}.
		\end{proof}
	\end{corollary}
	The numerical validation of the above corollary, showcasing the convergence of $\mathfrak{L}^{(\gamma)}_n(z,c;q)$ to $\mathcal{L}_n(z,\gamma,c)$ as $q$ tends to $1$, is presented in Table~\ref{zeros of generalized q-laguerre polynomials as q approaches to 1}. The table provides evidence of the convergence of zeros for the polynomials $\mathfrak{L}_6^{(0.1)}(z,1;q)$ as $q$ approaches $1$. These zeros closely align with the zeros of the polynomials $\mathcal{L}_6(z,\gamma=0.1,c=1)$. 
	
	\begin{table}[ht]
		\begin{center}
			\begin{tabular}{|c c c|c|}
				\hline\hline
				\multicolumn{3}{|c|}{Zeros of $\mathfrak{L}_6^{(0.1)}(z,1;q)$}&Zeros of\\
				$q=0.9999$&$q=0.99999$&$q=0.999990001$&$\mathcal{L}_6(z,\gamma=0.1,c=1)$\\
				\hline\hline
				0.606395&0.60605&0.606169&0.606014\\
				\hline
				1.851&1.84994&1.84993&1.84977\\
				\hline
				3.99236&3.98959&3.98979&3.98969\\
				\hline
				7.05592&7.04615&7.04316&7.0438\\
				\hline
				11.3481&11.3815&11.3962&11.387\\
				\hline
				17.9604&17.8383&17.8026&17.8237\\
				\hline
			\end{tabular}
		\end{center}
		\caption{Zeros of $\mathfrak{L}_6^{(\gamma=0.1)}(z,c=1;q)$ when $q \to 1$ converges to $\mathcal{L}_6(z,\gamma=0.1,c=1)$}
		\label{zeros of generalized q-laguerre polynomials as q approaches to 1}
	\end{table}
	Indeed, as indicated in \eqref{definition of Sobolev orthogonal polynomials}, when the inner product involves derivatives of polynomials, the polynomials are identified as Sobolev orthogonal polynomials.  In the $q$-analogue, this involves substituting derivatives with the $q$-difference operator. 
	
	\begin{definition}\label{defintion of sobolev orthogonality}
		Suppose that ${p_n(z;q)}_{n=0}^{\infty}$ constitutes a set of $q$-orthogonal polynomials associated with the positive measure $d_q\mu(z)$ over the interval $[a,b]$. Then, the set $\{y_n(z;q)\}_{n=0}^{\infty}$ becomes $q$-Sobolev orthogonal if it satisfies a $q$-difference equation of the structure given in \eqref{q-difference equation with required condition}.
	\end{definition}
	Hence, based on Theorems~\ref{relation between OP and Sobolev type OP for little q-Jacobi} and \ref{relation between OP and Sobolev type OP for q-laguerre}, we can affirm that the polynomials $\mathfrak{P}^{(\gamma,\xi)}_n(z,c;q)$ and $\mathfrak{L}^{(\gamma)}_n(z,c;q)$ satisfy the conditions for being $q$-Sobolev orthogonal polynomials.
	
	In addition to the mentioned $q$-polynomials, several other families of $q$-polynomials satisfy the hypergeometric $q$-difference equation \eqref{q-difference equation with required condition}. For instance, we have the generalized $q$-Bessel polynomials $\mathfrak{B}^{(\xi)}_n(z,c;q):={_3\phi_2}(q^{-n},-q^{\xi+n},q;0,q^{c+1};q,qz)$, the extended little $q$-Laguerre polynomials $\mathfrak{C}^{(\gamma)}_n(z,c;q):={_3\phi_2}(q^{-n},0,q;q^{\gamma+1},q^{c+1};q,qz)$, and the generalized Stieltjes-Wigert polynomials $\mathcal{S}_n(z,c;q):={_2\phi_2}(q^{-n},q;0,q^{c+1};q,-q^{n+1}z)$, which are derived from $q$-Bessel, little $q$-Laguerre, and Stieltjes-Wigert polynomials, respectively (see \cite{2010_Koekoek Lesky Swarttow_Hypergeometric orthogonal polynomials}). It is worth noting that these polynomials are also special cases of $\mathfrak{P}^{(\gamma,\xi)}_n(z,c;q)$ and $\mathfrak{L}^{(\gamma)}_n(z,c;q)$.
	
	These polynomials possess integral representations and satisfy mixed recurrence relations. We can obtain these relations by taking the limits of $\mathfrak{P}^{(\gamma,\xi)}_n(z,c;q)$ and $\mathfrak{L}^{(\gamma)}_n(z,c;q)$ under specific conditions as presented below.
	\begin{equation*}
		\begin{aligned}
			&\mathfrak{B}^{(\xi)}_n(z,c;q)=\lim_{\gamma\to \infty}\{\mathfrak{P}^{(\gamma,\xi)}_n(z,c;q):q^{\xi}\to-q^{-\gamma-1+\xi}\},\\
			&\mathfrak{C}^{(\gamma)}_n(z,c;q)=\lim_{\xi\to \infty}\mathfrak{P}^{(\gamma,\xi)}_n(z,c;q),\\
			&\mathcal{S}_n(z,c;q)=\lim_{\gamma\to \infty}\mathfrak{L}^{(\gamma)}_n(q^{-\gamma}z,c;q).
		\end{aligned}
	\end{equation*}
	With these notations, the following theorems can be immediate and the proof of $\mathfrak{B}^{(\xi)}_n(z,c;q),\\\mathfrak{C}^{(\gamma)}_n(z,c;q)$ follows from Theorem~\ref{relation between OP and Sobolev type OP for little q-Jacobi} and $\mathcal{S}_n(z,c;q)$ follows from Theorem~\ref{relation between OP and Sobolev type OP for little q-laguerre}. Hence, we omit the proofs and provide the statements only.
	\begin{theorem}{\label{relation between OP and Sobolev type OP for q-bessel}}
		For arbitrary $c\in\mathbb{N}$, the polynomials $y_n(z;q)=\mathfrak{B}^{(\xi)}_n(z,c;q)$, satisfy the orthogonality relation \eqref{definition of sobolev orthogonality-discrete} with \;$t=c$;\;$q^{\xi}>0$,\; $d\mu_q(z)=\sum_{k=0}^{\infty}\frac{q^{\xi k}}{(q;q)_k}q^{k(k+1)/2}\delta(z-q^k),$
		\vspace{-0.3em}
		\begin{align*}
			A_{n,q}=(q;q)_n(-q^{\xi+n};q)_{\infty}\frac{q^{\xi n}q^{n(n+1)/2}}{(1+q^{\xi+2n})}\bigg(\frac{(q;q)_c}{(1-q)^{c}}\bigg)^2, \; n\in \mathbb{N}_0.
		\end{align*}
		Furthermore, the polynomial $y_n(z;q)$ adheres to the $q$-difference equation \eqref{q-difference equation of Sobolev OP} with $\gamma'=2,\;\tilde{l}_2(z;q)=q^{\xi}z^2,\;\tilde{l}_1(z;q)=\frac{1+q^{\xi+1}}{q(q-1)}z+\frac{1}{q^2(1-q)},\;\tilde{l}_0(z;q)=0\;and\;\lambda_{n,q}=-\frac{[n]_q(1+q^{n+\xi})}{q^{n+1}(1-q)}.$
	\end{theorem}
	\begin{theorem}{\label{relation between OP and Sobolev type OP for little q-laguerre}}
		For arbitrary $c\in\mathbb{N}$, the polynomials $y_n(z;q)=\mathfrak{C}_n^{(\gamma)}(z,c;q)$, satisfy the orthogonality relation \eqref{definition of sobolev orthogonality-discrete} with \, $t=c;\;0<q^{\gamma+1}<1$, $d\mu_q(x)=\sum_{k=0}^{\infty}\frac{(q^{\gamma+1})^k}{(q;q)_k}\delta(z-q^k)$,
		\vspace{-0.2em}
		\begin{align*}
			A_{n,q}=\frac{(q^{\gamma+1})^n}{(q^{\gamma+1};q)_{\infty}}\frac{(q;q)_n}{(q^{\gamma+1};q)_{n}}\bigg(\frac{(q;q)_c}{(1-q)^{c}}\bigg)^2, \; n\in \mathbb{N}_0.
		\end{align*}
		Additionally, the polynomial $y_n(z;q)$ exhibits the q-difference equation \eqref{q-difference equation of Sobolev OP} with $\gamma'=2,\;\\\tilde{l}_2(z;q)=q^{\gamma}z,\;\tilde{l}_1(z;q)=\frac{1}{q-1}z+\frac{[\gamma+1]_q}{q},\;\tilde{l}_0(z;q)=0\;and\;\lambda_{n,q}=-\frac{[n]_q}{q^n(1-q)}.$
	\end{theorem}
	\begin{theorem}{\label{relation between OP and Sobolev type OP for Stieltjes-wigert}}
		For arbitrary $c\in\mathbb{N}$, the polynomials $y_n(z;q)=\mathcal{S}_n(z,c;q)$, satisfy the orthogonality relation \eqref{definition of sobolev orthogonality} with $t=c$;\; $K=[0,\infty)$, $w_q(z)=\frac{1}{(-z,-qz^{-1};q)_{\infty}}$,
		\vspace{-0.2em}
		\begin{align*}
			A_{n,q}=-\frac{lnq}{q^n}\frac{(q;q)_{\infty}}{(q;q)_{n}}\bigg(\frac{(q;q)_c}{(1-q)^{c}}\bigg)^2, \; n\in \mathbb{N}_0.
		\end{align*}
		Moreover, the polynomial $y_n(z;q)$ solves the q-difference equation \eqref{q-difference equation of Sobolev OP} with $\gamma'=2,\;\\\tilde{l}_2(z;q)=z^2,\;\tilde{l}_1(z;q)=\frac{1}{q-1}z+\frac{1}{q^2(1-q)},\;\tilde{l}_0(z;q)=0\;and\;\lambda_{n,q}=-\frac{[n]_q}{q(1-q)}.$
	\end{theorem}
	\begin{remark}
		From theorems~\ref{relation between OP and Sobolev type OP for q-bessel}, \ref{relation between OP and Sobolev type OP for little q-laguerre}, and \ref{relation between OP and Sobolev type OP for Stieltjes-wigert}, we can conclude that $\mathfrak{B}^{(\xi)}_n(z,c;q)$, $\mathfrak{C}^{(\gamma)}_n(z,c;q)$, and $\mathcal{S}_n(z,c;q)$ are $q$-Sobolev orthogonal polynomials.
	\end{remark}
	
	In the following discussion, we illustrate the categorization of $q$-Sobolev orthogonal polynomials into the class of classical orthogonal polynomials. The properties of classical orthogonal polynomials are elaborated in \cite{1194_Paco_Classical orthogonal polynomials}. However, in our current context of $q$-orthogonal polynomials, their characteristics are provided in \cite{2011_vinet_missing family}. Polynomials are classified as classical $q$-orthogonal polynomials if they are solutions of the TTRR and eigenvalue representation, as delineated in \cite{2011_vinet_missing family}. In particular, regarding little $q$-Jacobi, $q$-Bessel, Stieltjes-Wigert, little $q$-Laguerre, and $q$-Laguerre polynomials, they adhere to the TTRR, as outlined in \cite{2010_Koekoek Lesky Swarttow_Hypergeometric orthogonal polynomials}, alongside the eigenvalue equation outlined in Theorems~\ref{relation between OP and Sobolev type OP for little q-Jacobi}, \ref{relation between OP and Sobolev type OP for q-laguerre}, \ref{relation between OP and Sobolev type OP for q-bessel}, \ref{relation between OP and Sobolev type OP for little q-laguerre}, and \ref{relation between OP and Sobolev type OP for Stieltjes-wigert}. Consequently, we can affirm that these polynomials fall into the $q$-orthogonal polynomials of the classical category. Therefore, by fulfilling the eigenvalue equation, the system of polynomials ${y_n(z;q)}_{n=0}^{\infty}$ is characterized as the $q$-Sobolev orthogonal polynomials of classical type.
	
	\section{Behavior of zeros} \label{behavior of zeros}
	In this section, our focus is on analyzing the zeros of $\mathfrak{P}^{(\gamma,\xi)}_n(z,c;q)$ and $\mathfrak{L}^{(\gamma)}_n(z,c;q)$ polynomials, where $\gamma, \xi > -1$. The behavior of zeros for $-2 < \gamma, \xi < -1$ has been studied in \cite{2018_Kerstin_Zeros of Jacobi polynomials}, revealing insights into the quasi-orthogonal behavior of these polynomials. It has practical applications in various domains such as numerical analysis, approximation theory, and solving $q$-difference equations. To investigate these zeros, we use the matrix associated with the extended recurrence relation \eqref{recurrence relation for genralized little q-jacobi polynomials}.
	{\tiny{
			\begin{align*}
				\begin{bmatrix}
					\mu_{q,3}(0)+x\mu_{q,6}(0)&\mu_{q,4}(0)&&&\bf{0} \\
					\mu_{q,2}(1)+x\mu_{q,6}(1)&\mu_{q,3}(1)+x\mu_{q,6}(1) &\mu_{q,4}(1)\\
					\mu_{q,2}(2)&\mu_{q,2}(2)+x\mu_{q,5}(2)&\mu_{q,3}(2)+x\mu_{q,6}(2)\\
					&\mu_{q,1}(3)&\mu_{q,2}(3)+x\mu_{q,5}(3)\\
					&\ddots&\ddots&\ddots&\ddots\\
					\bf{0}
				\end{bmatrix}
				\begin{bmatrix}
					\mathfrak{P}_0^{(\gamma,\xi)}(z,c;q)\\
					\mathfrak{P}_1^{(\gamma,\xi)}(z,c;q)\\
					\mathfrak{P}_2^{(\gamma,\xi)}(z,c;q)\\
					\vdots\\
					\mathfrak{P}_{n-1}^{(\gamma,\xi)}(z,c;q)
				\end{bmatrix}
				+\mu_{4,q}(n-1)\begin{bmatrix}
					0\\
					0\\
					0\\
					\vdots\\
					\mathfrak{P}_n^{(\gamma,\xi)}(z,c;q)
				\end{bmatrix}
				=
				\bf{0}
	\end{align*}}}
	\begin{theorem}\label{Generalized eigenvalue problem}
		For any $n\ge2$,
		\begin{align*}
			{\bf A}_nX(x)=x {\bf B}_nX(x)-\mu_{4,q}(n-1)\mathfrak{P}_n^{(\gamma,\xi)}(z,c;q){\bf e}_n, 
		\end{align*}
		where 
		\begin{align*}
			&{\bf A}_n=\begin{bmatrix}
				\mu_{q,3}(0)&\mu_{q,4}(0)&&&&&\bf{0}\\
				\mu_{q,2}(1)&\mu_{q,3}(1)&\mu_{q,4}(1)\\
				\mu_{q,1}(2)&\mu_{q,2}(2)&\mu_{q,3}(2)&\mu_{q,4}(2)\\
				&\mu_{q,1}(3)&\mu_{q,2}(3)&\mu_{q,3}(3)\\
				&&\mu_{q,1}(4)\\
				&\ddots&\ddots&\ddots&\ddots\\
				\bf{0}
			\end{bmatrix}_{n\times n},
		\end{align*}
		\begin{align*}
			&{\bf B}_n=\begin{bmatrix}
				-\mu_{q,6}(0)&&&&&\bf{0}\\
				-\mu_{q,5}(1)&-\mu_{q,6}(1)\\
				&-\mu_{q,5}(2)&-\mu_{q,6}(2)\\
				&&-\mu_{q,5}(3)\\
				&&&\ddots&\ddots\\
				\bf{0}
			\end{bmatrix}_{n\times n},
		\end{align*}
		$ X(x)=(\mathfrak{P}_0^{(\gamma,\xi)}(z,c;q),\mathfrak{P}_1^{(\gamma,\xi)}(z,c;q),...,\mathfrak{P}_{n-1}^{(\gamma,\xi)}(z,c;q))^T$ and ${\bf e}_n=(0,0,...,1)^T$. Furthermore, the zeros of $\mathfrak{P}^{(\gamma,\xi)}_n(z,c;q)$ correspond to the eigenvalues of the generalized eigenvalue problem expressed as
		\begin{align*}
			{\bf A}_nX(z)=x {\bf B}_nX(z).
		\end{align*}
	\end{theorem}
	
	Hence, we have two methods for finding the zeros: we can either solve the hypergeometric series associated with the polynomials $\mathfrak{L}^{(\gamma)}_n(z,c;q)$ and $\mathfrak{P}^{(\gamma,\xi)}_n(z,c;q)$, or we can determine the eigenvalues through generalized eigenvalue problem.
	
	The generalized eigenvalue problem (see \cite{2021_Behera_generalized inverse eigenvalue problem,2019_Sri Ranga_ generalized eigenvalue problem}) provides a means to examine the properties of measures and orthogonal polynomials, as well as to study the relationships between them and the eigenvalues associated with specific linear operators.

	A detailed computation shows that $\mathfrak{P}^{(\gamma,\gamma)}_2(z,\gamma;q)$ and $\mathfrak{L}^{(\gamma)}_2(z,\gamma;q)$ have complex zeros for extremely large value of $\gamma$. We omit the details here to avoid lengthy expressions. In the context of a positive-definite linear functional, the roots (zeros) of orthogonal polynomials are simple, real, and typically situated within the specified interval of orthogonality. However, in certain scenarios where specific combinations of $\gamma$, $\xi$, and $c$ are employed, the resulting polynomials $\mathfrak{P}^{(\gamma,\gamma)}_2(z,\gamma;q)$ and $\mathfrak{L}^{(\gamma)}_2(z,\gamma;q)$ exhibit complex zeros that fall beyond this interval, clearly indicating a departure from orthogonality.
	
	In the subsequent analysis, we explore how the zeros of classical-type $q$-Sobolev orthogonal polynomials behave as all parameters approach the value of 1. The figures presented below provide a clear visual depiction that the vicinity of $\gamma=1$, $\xi=1$, and $c=1$ is a critical point influencing the behavior of zeros. Notably, the behavior of zeros transforms as the parameters approach the value of 1.
	
	\begin{figure}[h!]
		\footnotesize
		\stackunder[5pt]{\includegraphics[scale=0.55]{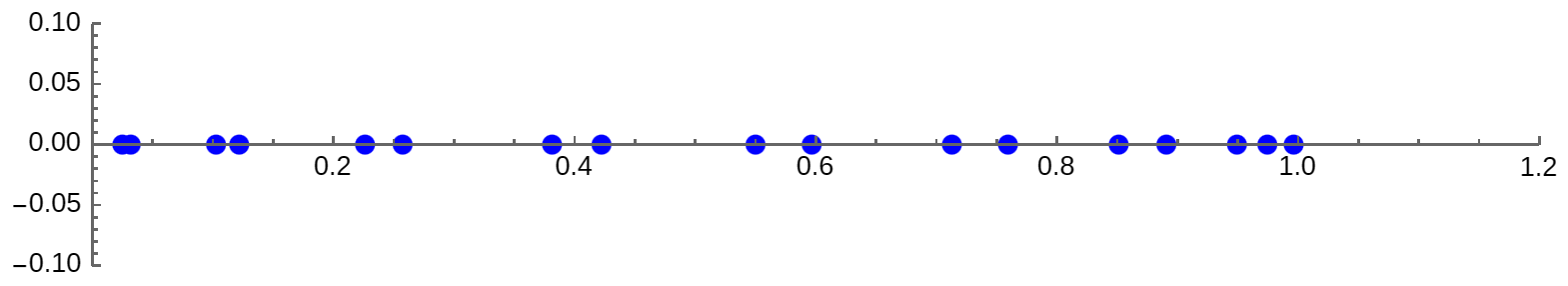}}{Zeros of $\mathfrak{P}^{(0.7,0.7)}_{17}(z,0.7;0.99998)$}
	\end{figure}
	\begin{figure}[h!]
		\footnotesize
		\stackunder[5pt]{\includegraphics[scale=0.7]{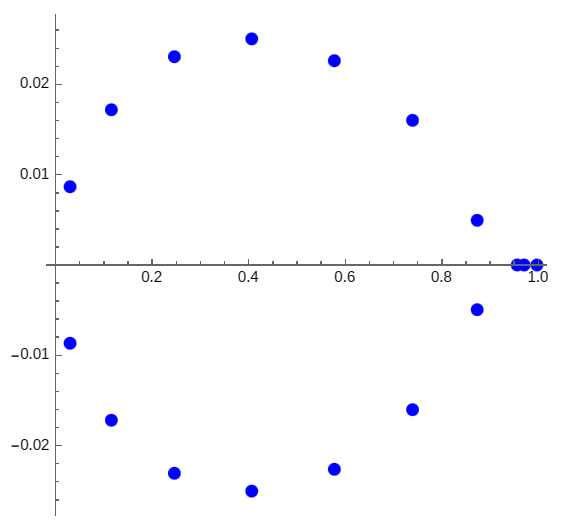}}{Zeros of $\mathfrak{P}^{(0.8,0.8)}_{17}(z,0.8;0.99998)$}
		\hspace{1cm}%
		\stackunder[5pt]{\includegraphics[scale=0.7]{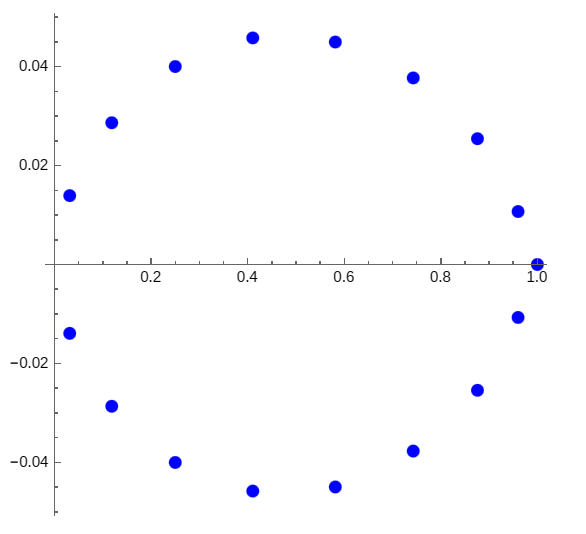}}{Zeros of $\mathfrak{P}^{(0.9,0.9)}_{17}(z,0.9;0.99998)$}
		\hspace{1cm}%
		\stackunder[5pt]{\includegraphics[scale=0.7]{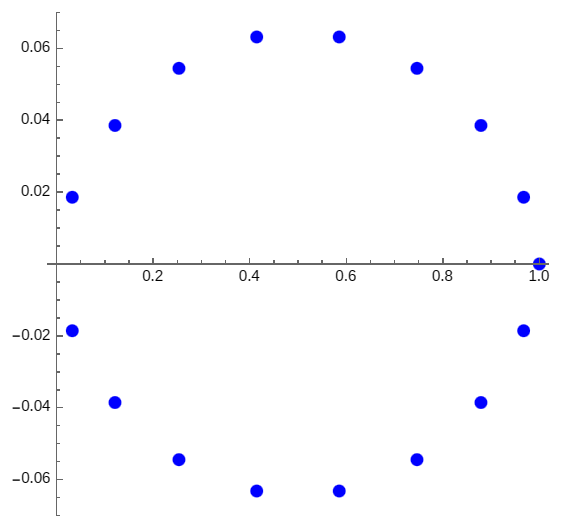}}{Zeros of $\mathfrak{P}^{(1,1)}_{17}(z,1;0.99998)$}
		\caption{Zeros of polynomials $\mathfrak{P}^{(\gamma,\xi)}_n(z,c;q)$ when $\gamma,\,\xi$ and $c\to1$.}
		\label{figure of generalized little q-jacobi polynomials}
	\end{figure}
	As evident in the zeros of $\mathfrak{P}_{17}^{(\gamma,\xi)}(z,c;0.99998)$, a discernible pattern emerges in the distribution of zeros as the parameter values are varied. All zeros align along the real axis when parameters are $\gamma=\xi=c = 0.7$. In contrast, for $\gamma=\xi=c = 0.8$, three zeros remain real while the rest become complez, as the figure~\ref{figure of generalized little q-jacobi polynomials} shows. Progressing to $\gamma=\xi=c = 0.9$, the arrangement of zeros becomes distorted, resembling a slightly oval shape in the right half. However, a remarkable transformation occurs when the parameters $\gamma$, $\xi$, and $c$ reach the value 1. At this critical point, the zeros of the polynomials converge and align to form a distinctive circular pattern. Notably, the zero assembly takes on a circular form for parameter values greater than one, albeit with a radius on the imaginary axis surpassing the previous value. As a result, $\gamma=\xi=c = 1$ serves as pivotal points governing the configuration of zeros. In summary, the figures provide insights into the behavior of zeros of polynomials and their critical points.
	\begin{figure}	
		\footnotesize
		\stackunder[5pt]{\includegraphics[scale=0.7]{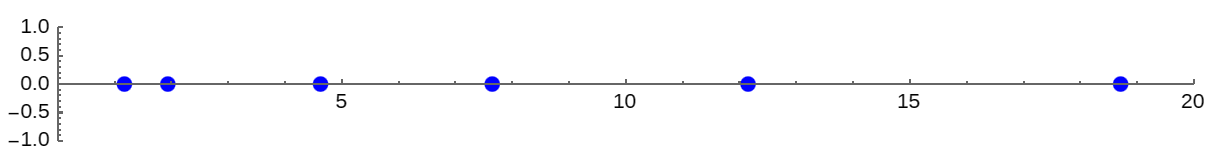}}{Zeros of $\mathfrak{L}^{(0.8)}_{6}(z,0.8;0.99999)$}
		\hspace{1cm}%
		\stackunder[5pt]{\includegraphics[scale=0.7]{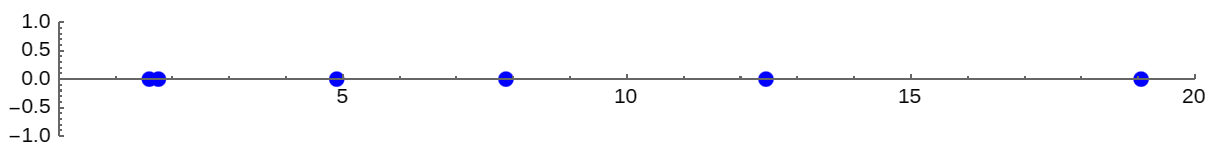}}{Zeros of $\mathfrak{L}^{(0.9)}_{6}(z,0.9;0.99999)$}
	\end{figure}
	\begin{figure}
		\footnotesize
		\stackunder[5pt]{\includegraphics[scale=0.7]{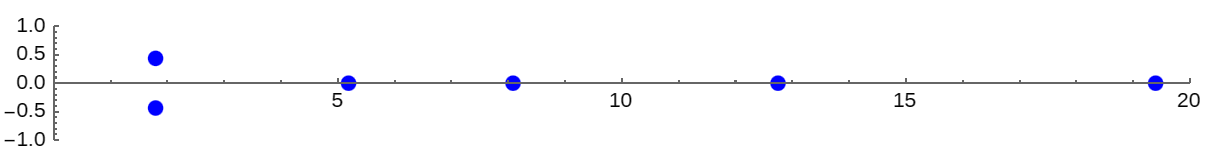}}{Zeros of $\mathfrak{L}^{(1)}_{6}(z,1;0.99999)$}
		\caption{Zeros of polynomials $\mathfrak{L}^{(\gamma)}_{6}(z,c;q)$ when $\gamma,\,c\to1$.}
	\end{figure}
	
	Also, in the case of $\mathfrak{L}^{(\gamma)}_{n}(z,c;q)$, the behavior of zeros undergo a distinct change as the parameters $\gamma$ and $\xi$ approach the value of 1 for the $\mathfrak{L}^{(\gamma)}_{6}(z,c;0.99999)$ polynomials. Specifically, there is a shift from real zeros for $\gamma=\xi=0.8$ and $\gamma=\xi=0.9$ to complex zeros at $\gamma=\xi=1$, and this complex behavior persists as the parameters continue to increase from 1. This observation emphasizes the alteration in the behavior of zeros as $\gamma$ and $c$ approach the value of 1.
	
	The interlacing behavior of zeros for various sequences corresponding to different parameter values of $q$-orthogonal polynomials has been studied and documented in \cite{2010_Jordaan_Mixed recurrence relations}. In our current discussion, we are focusing on the interlacing properties of zeros for the polynomials $\mathfrak{L}^{(\gamma)}_n(z,1;q)$ and $\mathcal{L}^{(\gamma)}_n(z;q)$.
	\begin{definition}\cite{1978_Chihara_orthogonal polynomials}
		A system of polynomials $\{P_n(z)\}_{n=0}^{\infty}$ with real roots $a_{n,1}\leq a_{n,2}\leq\ldots\leq a_{n,n}$ is said to have interlacing roots if, for all $n\geq 1$, the following holds
		\begin{align*}
			a_{n+1,1}\leq a_{n,1}\leq a_{n+1,2}\leq\ldots\leq a_{n,n}\leq a_{n+1,n+1},\quad n=0,1,\ldots.
		\end{align*}
	\end{definition}
	\subsection{Behavior of zeros of  $\mathfrak{L}_n^{(\gamma)}(z,c=1;q)$ for $-1<\gamma<1$}
	Table~\ref{zeros of generalized q-laguerre polynomials} lists the zeros of polynomials $\mathfrak{L}_6^{(0.5)}(z,1;0.9)$ and $\mathfrak{L}_7^{(0.5)}(z,1;0.9)$. In Figure ~\ref{fig:Generalizedq-Laguerreforq=0.9,alpha=0.5}, it is evident that these zeros exhibit the interlacing property. Thus, interlacing of zeros signifies that these polynomials demonstrate the separation theorem \cite{2007_jordaan_Separation theorems}, a notable characteristic found in specific families of orthogonal polynomials.
	\begin{figure}[h!]
		\includegraphics[scale=0.6]{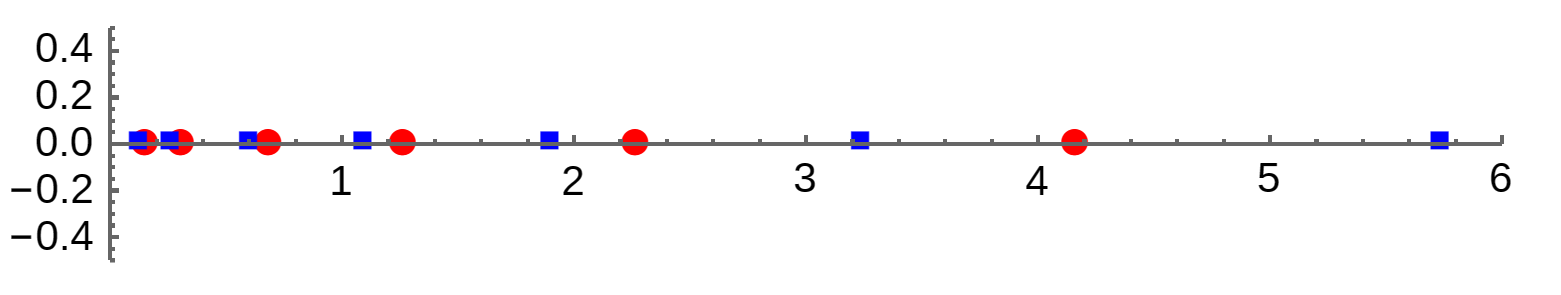}
		\caption{Zeros of $\mathfrak{L}_7^{(0.5)}(z,1;0.9)$ (blue squares) and $\mathfrak{L}_6^{(0.5)}(z,1;0.9)$ (red circles).}
		\label{fig:Generalizedq-Laguerreforq=0.9,alpha=0.5}
	\end{figure}
	\begin{table}[ht]
		\begin{center}
			\begin{tabular}{|c|c|}
				\hline
				Zeros of $\mathfrak{L}_6^{(0.5)}(z,1;0.9)$ &Zeros of $\mathfrak{L}_7^{(0.5)}(z,1;0.9)$\\
				\hline
				0.127251&0.117823\\
				\hline
				0.282878&0.25455\\
				\hline
				0.659823&0.59272\\
				\hline
				1.23947&1.08522\\
				\hline
				2.2416&1.89115\\
				\hline
				4.13778&3.23078\\
				\hline
				-&5.72914\\
				\hline
			\end{tabular}
		\end{center}
		\caption{}
		\label{zeros of generalized q-laguerre polynomials}
	\end{table}
	\begin{lemma}
		The polynomials $\mathfrak{L}^{(\gamma)}_n(z,1;q)$ for $n\in \mathbb{N}$ and $\gamma>-1$ exhibit the following relation
		\begin{equation}\label{relation between generalized q-laguerre polynomials}
			\begin{aligned}
				(1+x)\mathfrak{L}^{(\gamma)}_n(qz,1;q)&=\frac{(1-q^{n+\gamma+1})(1-q^{n+2})}{q^{n+\gamma+1}(q^{n+1}-1)}\mathfrak{L}^{(\gamma)}_{n+1}(z,1;q)\\
				&+\frac{q^{2n+\gamma+1}+1-q^{2n+\gamma+2}}{q^{n+\gamma+1}}\mathfrak{L}^{(\gamma)}_n(z,1;q).
			\end{aligned}
		\end{equation}
		\begin{proof}
			This result is achieved by comparing coefficients of $x^n$ on both sides of \eqref{relation between generalized q-laguerre polynomials}.
		\end{proof}
	\end{lemma}
	
	\begin{theorem}{\label{interlacing between sobolev q-laguerre and q-laguerre}}
		Let $x_{n,j}$ and $y_{n,j}$ denote the zeros of the polynomials $\mathfrak{L}^{(\gamma)}_n(z,1;q)$ and $\mathcal{L}^{(\gamma)}_n(z;q)$, respectively, for $1 \leq j \leq n$. Then, in the case of $-1<\gamma<1$, the arrangement of these zeros is as follows
		\begin{align*}
			0\leq y_{n,1}\leq x_{n,1}\leq y_{n,2}\leq...\leq y_{n,n}\leq x_{n,n}.
		\end{align*} 
		\begin{proof}
			The following expression from Theorem~\ref{relation between OP and Sobolev type OP for q-laguerre} is known.
			\begin{align*}
				\mathcal{D}_q^{c}(z^{c}\mathfrak{L}_n^{(\gamma)}(z,c;q))=\frac{(q;q)_{c}}{(1-q)^{c}}\mathcal{L}_n^{(\gamma)}(z;q).
			\end{align*}
			For $c=1$, this becomes
			\begin{align*}
				q\mathfrak{L}^{(\gamma)}_n(qz,1;q)-\mathfrak{L}^{(\gamma)}_n(z,1;q)=(q-1)\mathcal{L}_n(z;\gamma,1).
			\end{align*}
			Using \eqref{relation between generalized q-laguerre polynomials}, we get
			\begin{equation}\label{relation between sobolev q-laguerre and q-laguerre}
				\begin{aligned}
					-\frac{(1-q)}{q}\mathcal{L}_n^{(\gamma)}(z;q)&=-\frac{(1-q^{n+\gamma+1})(1-q^{n+2})}{q^{n+\gamma+1}(1-q^{n+1})(1+z)}\mathfrak{L}^{(\gamma)}_{n+1}(z,1;q)+\\
					&+\frac{q^{2n+\gamma+1}+1-q^{2n+\gamma+2}}{q^{n+\gamma+1}(1+z)}\mathfrak{L}^{(\gamma)}_n(z,1;q)-\frac{1}{q}\mathfrak{L}^{(\gamma)}_n(z,1;q).
				\end{aligned}
			\end{equation}
			By evaluating the expression \eqref{relation between sobolev q-laguerre and q-laguerre} at successive zeros $x_{n,j}$ and $x_{n,j+1}$ of $\mathfrak{L}^{(\gamma)}_n(z,1;q)$, we obtain
			\begin{equation*}
				\begin{aligned}
					\frac{(1-q)^2}{q^2}(1+x_{n,j})(1+x_{n,j+1})\mathcal{L}_n^{(\gamma)}(x_{n,j};q)\mathcal{L}_n^{(\gamma)}(x_{n,j+1};q)&=\Bigg(\frac{(1-q^{n+\gamma+1})(1-q^{n+2})}{q^{n+\gamma+1}(1-q^{n+1})(1+z)}\Bigg)^2\\
					&\mathfrak{L}^{(\gamma)}_{n+1}(x_{n,j},1;q)\mathfrak{L}^{(\gamma)}_{n+1}(x_{n,j+1},1;q).
				\end{aligned}
			\end{equation*}
			Because the zeros of polynomials $\mathfrak{L}^{(\gamma)}_{n}(z,1;q)$ and $\mathfrak{L}^{(\gamma)}_{n+1}(z,1;q)$ for $-1<\gamma<1$ interlace, as illustrated in the Figure~\ref{fig:Generalizedq-Laguerreforq=0.9,alpha=0.5}, the right side expression of above relation is negative and therefore,
			\begin{align*}
				x_{n,j}\leq y_{n,j+1}\leq x_{n,j+1}\quad for\quad j=1,2,...,n-1.
			\end{align*}
			Hence, we obtained the desired result.
		\end{proof}
	\end{theorem}
	\begin{figure}[h!]
		\includegraphics[scale=0.65]{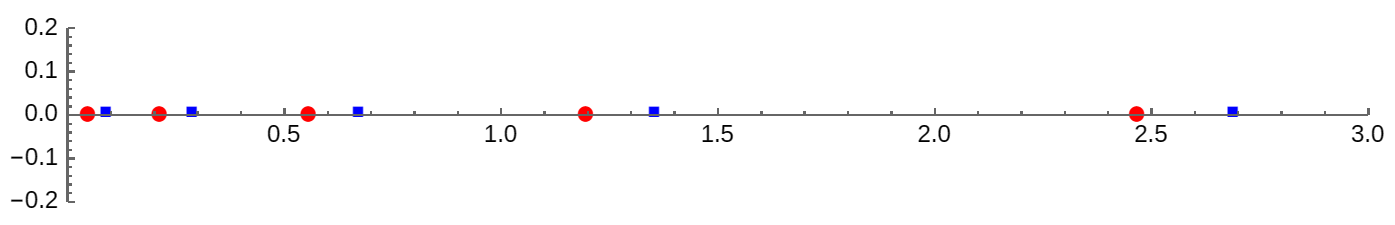}
		\caption{Zeros of $\mathfrak{L}^{(0.1)}_5(z,1;0.9)$ (blue squares) and $\mathcal{L}^{(0.1)}_5(z;0.9)$ (red circles).}
		\label{fig:Interlacing between sobolev q-laguerre and q-laguerre}
	\end{figure}
	
	\begin{minipage}[c]{0.32\textwidth}
		Table~\ref{tab:interlacing betweeen sobolev q-laguerre and q-laguerre} displays the zeros of the polynomials $\mathcal{L}^{(0.1)}_5(z;0.9)$ and $\mathfrak{L}^{(0.1)}_5(z,1;0.9)$. These zeros are also observed to be interlacing, as depicted in Figure~\ref{fig:Interlacing between sobolev q-laguerre and q-laguerre}. This numerical observation serves to validate the theorem outlined above.
	\end{minipage}
	\hfil
	\begin{minipage}[c]{0.4\textwidth}
		\centering
		\begin{tabular}{|c|c|}
			\hline
			Zeros of $\mathcal{L}^{(0.1)}_5(z;0.9)$ &Zeros of $\mathfrak{L}^{(0.1)}_5(z,1;0.9)$\\
			\hline
			0.039398&0.0874766\\
			\hline
			0.205233&0.285522\\
			\hline
			0.548521&0.669265\\
			\hline
			1.18799&1.35207\\
			\hline
			2.45937&2.68648\\
			\hline
		\end{tabular}
		\captionof{table}{}
		\label{tab:interlacing betweeen sobolev q-laguerre and q-laguerre}
	\end{minipage}
	\section{Concluding remarks}
	It is worth noting that having familiarity with the basic hypergeometric functions associated with these polynomials enable us to derive a third-order $q$-difference equation, as demonstrated in propositions~\ref{theorem of higher order q-difference equation for little q-Jacobi polynomials} and \ref{theorem of higher order q-difference equation for q-laguerre polynomials}. The order of these $q$-difference equations are in contrast to Theorems~\ref{relation between OP and Sobolev type OP for little q-Jacobi} and \ref{relation between OP and Sobolev type OP for laguerre}, where the $q$-difference equation is inherited directly from the second order $q$-difference equation of the $q$-Laguerre and little $q$-Jacobi polynomials. In the latter case, the order of the $q$-difference equation depends on the parameter $c$ value. We can infer the order from \eqref{q-difference equation of Sobolev OP}. These arguments lead to the following questions.
	\begin{problem}
		Can we determine the orthogonality properties of these $q$-Sobolev orthogonal polynomials using the information presented in Theorem~\ref{Generalized eigenvalue problem}?
	\end{problem}
	\begin{problem}
		Can a TTRR and five-term recurrence relation be derived for these $q$-Sobolev orthogonal polynomials under specific conditions?
	\end{problem}
	\section{Proof of Theorem~\ref{theorem for recurrence relation for generalized little q-jacobi polynomials}}\label{proof of theorem}
	Consider a positive integer, $\sigma$. We generalize the polynomials as follows:
	\begin{align*}
		&\textbf{L}_n^{(\gamma)}(z;q)=\textbf{L}_n^{(\gamma)}(z,c_1,c_2,...,c_{\sigma};q)\\
		&={_{\sigma+1}\phi _{\sigma+1}}(q^{-n},q,q,...,q;q^{\gamma+1},q^{c_1+1},q^{c_2+1},...,q^{c_{\sigma}+1};q;-q^{n+\gamma+1}z),\\ \label{generalized q-laguerre polynomials of order sigma}
		&\textbf{P}_n^{(\gamma,\xi)}(z;q)=\textbf{P}_n^{(\gamma,\xi)}(z,c_1,c_2,...,c_{\sigma};q)\\
		&={_{\sigma+2}\phi _{\sigma+1}}(q^{-n},q^{n+\gamma+\xi+1},q,...,q;q^{\gamma+1},q^{c_1+1},q^{c_2+1},...,q^{c_{\sigma}+1};q;qz), \\
		&\; \gamma, \xi>-1,\;c_1,\;...,c_{\sigma}>0,\; n \in \mathbb{N}_0.
	\end{align*}
	Fix an arbitrary integer $n \geq 2$. We can write
	\begin{equation*}
		\begin{aligned}
			\textbf{P}_{n+1}^{(\gamma,\xi)}(z;q)&=\sum_{i=0}^{n+1}\frac{{{(q^{-(n+1)};q)}_i}{{(q^{n+\gamma+\xi+2};q)}_i}{{(q;q)}_i}...{{(q;q)}_i}}{{{(q^{\gamma+1};q)}_i}{{(q^{c_1+1};q)}_i}...{{(q^{c_{\sigma}+1};q)}_i}}\frac{(qz)^i}{{{(q;q)}_i}}\\
			&=\sum_{i=0}^{n+1}\frac{{{(q^{-(n+1)};q)}_i}{{(q^{n+\gamma+\xi+2};q)}_i}{{(q;q)}_i}^{\sigma-1}}{{{(q^{\gamma+1};q)}_i}{{(q^{c_1+1};q)}_i}...{{(q^{c_{\sigma}+1};q)}_i}}(qz)^i\\
			&=\sum_{i=0}^{n+1}\tau_{i,q}(n),
		\end{aligned}
	\end{equation*}
	where
	\begin{align*}
		\tau_{i,q}(n)=\frac{{{(q^{-(n+1)};q)}_i}{{(q^{n+\gamma+\xi+2};q)}_i}{{(q;q)}_i}^{\sigma-1}}{{{(q^{\gamma+1};q)}_i}{{(q^{c_1+1};q)}_i}...{{(q^{c_{\sigma}+1};q)}_i}}(qz)^i.
	\end{align*}
	
	We must describe $\textbf{P}_n^{(\gamma,\xi)}(z;q),\; \textbf{P}_{n-1}^{(\gamma,\xi)}(z;q),\; \textbf{P}_{n-2}^{(\gamma,\xi)}(z;q),\; x\textbf{P}_n^{(\gamma,\xi)}(z;q) \; and \; x\textbf{P}_{n-1}^{(\gamma,\xi)}(z;q)$ in terms of $\tau_{i,q}(n)$. For any arbitrarily chosen $b \in \mathbb{R}$, we define the following:
	\begin{align*}
		&{{[b;q]}_0}=0,\,{{[b;q]}_1}=(1-b)\quad and\\ &{{[b;q]}_i}=(1-b)(1-bq^{-1})(1-bq^{-2})...(1-bq^{-i+1}),\;i=2,3,....
	\end{align*}
	\vspace{-1.21em}
	Using
	\begin{align*}
		{{(q^{-n};q)}_i}&={{(q^{-(n+1)};q)}_i}\frac{(1-q^{-n+i-1})}{(1-q^{-n-1})}\quad and\\
		{{(q^{n+\gamma+\xi+1};q)}_i}&={{(q^{n+\gamma+\xi+2};q)}_i}\frac{(1-q^{n+\gamma+\xi+1})}{(1-q^{n+\gamma+\xi+i+1})}, \; i \in \mathbb{N}_0,
	\end{align*}
	we obtain relation
	\begin{align}
		\textbf{P}_{n-t}^{(\gamma,\xi)}(z;q)=\sum_{i=0}^{\infty}\tau_{i,q}(n)q^{(t+1)i}\frac{{{[q^{n-i+1};q]}_{t+1}}{{[q^{n+\gamma+\xi+1};q]
				}_{t+1}}}{{{[q^{n+1};q]}_{t+1}}{{[q^{n+\gamma+\xi+i+1};q]}_{t+1}}},\quad t=0,1,2. \label{comparable expression 1}
	\end{align}
	Using
	\begin{align*}
		{{(a;q)}_{i-1}}=\frac{1}{(1-aq^{i-1})}{{(a;q)}_i}, \quad i \in \mathbb{N},
	\end{align*}
	we can deduce that
	\begin{equation}
		\begin{aligned}
			x\textbf{P}_{n-t}^{(\gamma,\xi)}(z;q)=\sum_{i=1}^{\infty}(-1)\tau_{i,q}(n)q^{n+(i-1)t}(1-q^i)\frac{{{[q^{n-i+1};q]}_{t}}{{[q^{n+\gamma+\xi+1};q]}_{t+1}}}{{{[q^{n+1};q]}_{t+1}}{{[q^{n+\gamma+\xi+i+1};q]}_{t+2}}}\\
			\times\frac{(1-q^{\gamma+i})(1-q^{c_1+i})...(1-q^{c_{\sigma}+i})}{(1-q^i)^{\sigma}},\quad t=0,1. \label{Comparable expression 2}
		\end{aligned}
	\end{equation}
	Consider the provided expression
	\begin{align*}
		Q(z;q)&=\mu_{q,1}(n)\textbf{P}_{n-2}^{(\gamma,\xi)}(z;q)+\mu_{q,2}(n)
		\textbf{P}_{n-1}^{(\gamma,\xi)}
		(z;q)+\mu_{q,3}(n)\textbf{P}_n^{(\gamma,\xi)}(z;q)+\mu_{q,4}(n)\\
		&\times\textbf{P}_{n+1}^{(\gamma,\xi)}(z;q)+\mu_{q,5}(n)x\textbf{P}_{n-1}
		^{(\gamma,\xi)}(z;q)+\mu_{q,6}(n)x\textbf{P}_{n}^{(\gamma,\xi)}(z;q),
	\end{align*}
	where $\mu_{q,i}(n)$,\;$1\leq i\leq 6$ denotes real free parameters. We choose these parameters to ensure $Q(z;q)=0$ for all $z$. $Q(z;q)$ takes the following form as a result of relations \eqref{comparable expression 1} and \eqref{Comparable expression 2}.
	\begin{align*}
		Q(z;q)&=\mu_{q,1}(n)\sum_{i=0}^{\infty}\tau_{i,q}(n)q^{3i}\frac{[q^{n-i+1};q]_3[q^{n+\gamma+\xi+1};q]_3}{[q^{n+1};q]_3[q^{n+\gamma+\xi+i+1};q]_3}+\mu_{q,2}(n)\sum_{i=0}^{\infty}\tau_{i,q}(n)q^{2i}\\
		&\times\frac{[q^{n-i+1};q]_2}{[q^{n+1};q]_2}\frac{[q^{n+\gamma+\xi+1};q]_2}{[q^{n+\gamma+\xi+i+1};q]_2}+\mu_{q,3}(n)\sum_{i=0}^{\infty}\tau_{i,q}(n)q^i\frac{[q^{n-i+1};q][q^{n+\gamma+\xi+1};q]}{[q^{n+1};q][q^{n+\gamma+\xi+i+1};q]}\\
		&+\mu_{q,4}(n)\sum_{i=0}^{\infty}\tau_{i,q}(n)+\mu_{q,5}(n)\sum_{i=1}^{\infty}(-1)\tau_{i,q}(n)q^{n+i-1}(1-q^i)\frac{[q^{n-i+1};q]_1}{[q^{n+1};q]_2}\\
		&\times\frac{[q^{n+\gamma+\xi+1};q]_2}{[q^{n+\gamma+\xi+i+1};q]_3}\frac{(1-q^{\gamma+i})(1-q^{c_1+i})...(1-q^{c_{\sigma}+i})}{(1-q^i)^{\sigma}}+\mu_{q,6}(n)\sum_{i=1}^{\infty}(-1)\tau_{i,q}(n)\\
		&\times q^{n}(1-q^i)\frac{[q^{n+\gamma+\xi+1};q]}{[q^{n+1};q][q^{n+\gamma+\xi+i+1};q]_2}\frac{(1-q^{\gamma+i})(1-q^{c_1+i})...(1-q^{c_{\sigma}+i})}{(1-q^i)^{\sigma}}.
	\end{align*}
	\\
	It can be simplified as 
	\begin{align*}
		Q(z;q)&=\tau_{0,q}(n)\big(\mu_{q,1}(n)+\mu_{q,2}(n)+\mu_{q,3}(n)+\mu_{q,4}(n)\big)+\sum_{i=1}^{\infty}\tau_{i,q}(n)A(\gamma,\xi,n,q;i),
	\end{align*}
	where
	\begin{align*}
		A(\gamma,\xi,n,q;i)&=\mu_{q,1}(n)q^{3i}\times\frac{[q^{n-i+1};q]_3[q^{n+\gamma+\xi+1};q]_3}{[q^{n+1};q]_3[q^{n+\gamma+\xi+i+1};q]_3}+\mu_{q,2}(n)q^{2i}\\
		&\times\frac{[q^{n-i+1};q]_2[q^{n+\gamma+\xi+1};q]_2}{[q^{n+1};q]_2[q^{n+\gamma+\xi+i+1};q]_2}+\mu_{q,3}(n)q^{i}\times\frac{[q^{n-i+1};q][q^{n+\gamma+\xi+1};q]}{[q^{n+1};q][q^{n+\gamma+\xi+i+1};q]}\\
		&+\mu_{q,4}(n)+\Big(\mu_{5,q}(n)q^{i-1}\frac{[q^{n-i+1};q](1-q^{n+\gamma+\xi})}{(1-q^n)(1-q^{n+\gamma+\xi+i-1})}+\mu_{q,6}(n)\Big)\\
		&\times(-1)\frac{[q^{n+\gamma+\xi+1};q]}{[q^{n+1};q][q^{n+\gamma+\xi+i+1};q]_2}\frac{(1-q^{\gamma+i})(1-q^{c_1+i})...(1-q^{c_{\sigma}+i})}{(1-q^i)^{\sigma}}\\
		&\times q^n(1-q^i).
	\end{align*}
	Observe that $\tau_{i,q}(n)=0$ for all $i\geq n + 2$. If
	\begin{align}
		\mu_{q,1}(n)+\mu_{q,2}(n)+\mu_{q,3}(n)+\mu_{q,4}(n)=0, \label{condition on first four mu functions}
	\end{align}
	and
	\begin{align}
		A(\gamma,\xi,n,q;i)=0\quad for \quad 1\leq i\leq n+1, \label{A(alpha,beta,n,q;i)=0}
	\end{align}
	then $Q(z;q)=0$. $A(\gamma,\xi,n,q;i)=\frac{A_i(n;q)}{B_i(n;q)}$, where $A_i(n;q)$ and $B_i(n;q)$ are polynomials of $q^i$, is a rational function of $q^i$. As a result, it is sufficient to check the equality $A_i(n;q)=0$ for a fixed number $M (M>degA_i(n;q))$ with distinct points $i$.
	
	In the following, we will assume that $\sigma=1, c_1=c>0$. Denote
	\begin{align*}
		\nu_{q,1}(n)&=\frac{[q^{n+\gamma+\xi+1};q]_3}{[q^{n+1};q]_3}\mu_{q,1}(n),& \nu_{q,2}(n)&=\frac{[q^{n+\gamma+\xi+1};q]_2}{[q^{n+1};q]_2}\mu_{q,2}(n),\\
		\nu_{q,3}(n)&=\frac{[q^{n+\gamma+\xi+1};q]}{[q^{n+1};q]}\mu_{q,3}(n),& \nu_{q,4}(n)&=\mu_{4,q}(n),\\
		\nu_{q,5}(n)&=(-1)\frac{[q^{n+\gamma+\xi+1};q]_2}{[q^{n+1};q]_2}\mu_{q,5}(n),& \nu_{q,6}(n)&=(-1)\frac{[q^{n+\gamma+\xi+1};q]}{[q^{n+1};q]}\mu_{q,6}(n).
	\end{align*}
	Now, we multiply equation \eqref{A(alpha,beta,n,q;i)=0} by $[q^{n+\gamma+\xi+i+1};q]_3$ to get 
	\begin{equation}
		\begin{aligned}
			&\nu_{q,1}(n)q^{3i}[q^{n-i+1};q]_3+\nu_{q,2}(n)q^{2i}(1-q^{n+\gamma+\xi+i-1})[q^{n-i+1};q]_2\\
			&+\nu_{q,3}(n)q^i[q^{n-i+1};q][q^{n+\gamma+\xi+i-1};q]_2+\nu_{q,4}(n;q)
			[q^{n+\gamma+\xi+i+1};q]_3\\
			&+\nu_{q,5}(n)q^{n+i-1}[q^{n-i+1};q](1-q^{\gamma+i})(1-q^{c+i})\\
			&+\nu_{q,6}(n)q^n(1-q^{n+\gamma+\xi+i-1})(1-q^{\gamma+i})(1-q^{c+i})=0, \quad 1\leq i\leq n+1. \label{Polynomial of degree less than 3 of i}
		\end{aligned}
	\end{equation}
	The expression on the left side represents a polynomial with a maximum degree of 3 in terms of $q^i$. Therefore, if \eqref{Polynomial of degree less than 3 of i} is valid for four distinct values of $i$, it will be valid for all complex values of $i$. It is important to choose these values carefully to keep the expressions as straightforward as possible. We select $i$ to be equal to $n-1$, $n+1$, $n$, and $0$. Using $i=n-2$ would result in a more complex expression. After substituting these chosen values and performing simplifications, including the utilization of \eqref{condition on first four mu functions}, we obtain
	\begin{align}
		&\nu_{q,4}(n)[q^{2n+\gamma+\xi+2};q]_2+\nu_{q,6}(n)(1-q^{\gamma+n+1})
		(1-q^{c+n+1})q^n=0,\notag\\
		&\nu_{q,3}(n)q^n[q;q][q^{2n+\gamma+\xi};q]_2+\nu_{q,4}(n)[q^{2n+\gamma+\xi+1};q]_3+\nu_{q,5}(n)q^{2n-1}[q;q](1-q^{\gamma+n})\notag\\
		&\times(1-q^{c+n})+\mu_{q,6}(n)(1-q^{2n+\gamma+\xi-1})(1-q^{\gamma+n})(1-q^{c+n})q^n=0.\notag\\
		\mbox{Also}\quad	&\nu_{q,2}(n)q^{2(n-1)}(1-q^{2n+\gamma+\xi-2})[q^2;q]_2+\nu_{q,3}(n)q^{n-1}[q^2;q][q^{2n+\gamma+\xi-1};q]_2\\
		&+\nu_{q,4}(n)[q^{2n+\gamma+\xi};q]_3+\nu_{q,5}(n)q^{2(n-1)}[q^2;q](1-q^{\gamma+n-1})(1-q^{c+n-1})\notag\\
		&+\nu_{q,6}(n)(1-q^{2n+\gamma+\xi-2})(1-q^{\gamma+n-1})(1-q^{c+n-1})q^n=0\notag
	\end{align}
	and
	\begin{align*}
		&\nu_{q,5}(n)q^{n-1}[q^{n+1};q](1-q^{\gamma})(1-q^{c})+\nu_{q,6}
		(n)(1-q^{n+\gamma+\xi-1})(1-q^{\gamma})(1-q^{c})q^n=0.
	\end{align*}
	We can deduce $\nu_{q,2}(n) ,\nu_{q,3}(n),\nu_{q,4}(n)$ and $\nu_{q,5}(n)$ in terms of $\nu_{q,6}(n)$ from the above expressions as
	\begin{align*}
		\nu_{q,4}(n)&=-\frac{q^n(1-q^{\gamma+n+1})(1-q^{c+n+1})}{[q^{2n+\gamma+\xi+2};q]_2}\nu_{q,6}(n),\\
		\nu_{q,5}(n)&=-\frac{q(1-q^{n+\gamma+\xi-1})}{[q^{n+1};q]}\nu_{q,6}(n)
		,\\
		\nu_{q,3}(n)&=\Bigg(\frac{(1-q^{c+n+1})(1-q^{\gamma+n+1})}{(1-q)
			(1-q^{2n+\gamma+\xi+2})}+\frac{q^n(1-q^{n+\gamma+\xi-1})(1-q^{\gamma+n})(1-q^{c+n})}{[q^{2n+\gamma+\xi};q]_2[q^{n+1};q]}\\
		&-\frac{(1-q^{\gamma+n})(1-q^{c+n})}{(1-q)(1-q^{2n+\gamma+\xi}
			)}\Bigg)\nu_{q,6}(n)\quad \mbox{and}\\
		\nu_{q,2}(n)&=\Bigg(-\frac{(1-q^{2n+\gamma+\xi-1})(1-q^{\gamma+n+1}
			)(1-q^{c+n+1})}{q^{n-1}(1-q)^{2}(1-q^{2n+\gamma+\xi+2})}-\frac{(1-q^{\gamma+n-1})(1-q^{c+n-1})}{q^{n-2}[q^2;q]_2}\\
		&+\frac{(1-q^{2n+\gamma+\xi-1})(1-q^{\gamma+n})(1-q^{c+n})}{q^{n-1}(1-q)^{2}(1-q^{2n+\gamma+\xi})}+\frac{[q^{2n+\gamma+\xi}
			;q]_2(1-q^{\gamma+n+1})(1-q^{c+n+1})}{q^{n-2}[q^2;q]_2[q^{2n+\gamma+\xi+2};q]_2}\\
		&+\frac{q(1-q^{\gamma+n-1})(1-q^{n+\gamma+\xi-1})(1-q^{c+n-1})
		}{(1-q)(1-q^{2n+\gamma+\xi-2})[q^{n+1};q]}\\
		&-\frac{q(1-q^{n+\gamma+\xi-1})(1-q^{\gamma+n})(1-q^{c+n})}{(1-q)(1-q^{2n+\gamma+\xi})[q^{n+1};q]}\Bigg)\nu_{q,6}(n).
	\end{align*}
	For ease of calculation, we now choose
	\begin{align*}
		\mu_{q,6}(n)=[q^{2n+\gamma+\xi+2};q]_2(1-q^{n+\gamma+\xi
		})(1-q^{n+1}).
	\end{align*}
	Then
	\begin{equation}\label{Coefficients of recurrence realtions}
		\begin{aligned}
			\mu_{q,3}(n)&=(-1)[q^{2n+\gamma+\xi+1};q][q^{n+\gamma+\xi};q]\Bigg(\frac{[q^{n+1};q](1-q^{\gamma+n+1})(1-q^{c+n+1})}{(1-q)}\\
			&+\frac{q^n(1-q^{n+\gamma+\xi-1})(1-q^{\gamma+n})(1-q^{c+n})(1-q^{2n+\gamma+\xi+2})}{[q^{2n+\gamma+\xi};q]_2}\\
			&-\frac{(1-q^{\gamma+n})(1-q^{c+n})(1-q^{2n+\gamma+\xi+2})[q^{n+1};q]}{(1-q)(1-q^{2n+\gamma+\xi})}\Bigg),\\
			\mu_{q,4}(n)&=q^n(1-q^{\gamma+n+1})(1-q^{c+n+1})[q^{n+\gamma+\xi+1};q]_2,\\
			\mu_{q,5}(n)&=(-1)q[q^n;q](1-q^{n+\gamma+\xi-1})[q^{2n+\gamma+\xi+2};q]_2,\\
			\mu_{q,2}(n)&=\frac{(1-q^{2n+\gamma+\xi-1})[q^{n+1};q]_2(1-q^{\gamma+n+1})(1-q^{c+n+1})(1-q^{2n+\gamma+\xi+1})}{q^{n-1}(1-q)^{2}}\\
			&+\frac{q(1-q^{n+\gamma+\xi-1})(1-q^{\gamma+n})(1-q^{c+n})(1-q^n)[q^{2n+\gamma+\xi+2};q]_2}{(1-q)(1-q^{2n+\gamma+\xi})}\\
			&-\frac{(1-q^{2n+\gamma+\xi-1})(1-q^{\gamma+n})(1-q^{c+n})[q^{n+1};q]_2[q^{2n+\gamma+\xi+2};q]_2}{q^{n-1}(1-q)^{2}(1-q^{2n+\gamma+\xi})}\\
			&-\frac{[q^{2n+\gamma+\xi};q]_2(1-q^{\gamma+n+1})(1-q^{c+n+1})[q^{n+1};q]_2}{q^{n-2}[q^2;q]_2}\\
			&-\frac{q(1-q^{n+\gamma+\xi-1})(1-q^{\gamma+n-1})(1-q^{c+n-1})(1-q^n)[q^{2n+\gamma+\xi+2};q]_2}{(1-q)(1-q^{2n+\gamma+\xi-2})}\\
			&+\frac{[q^{2n+\gamma+\xi+2};q]_2(1-q^{\gamma+n-1})(1-q^{c+n-1})[q^{n+1};q]_2}{q^{n-2}[q^2;q]_2}\\
			\mbox{and}\quad
			\mu_{q,1}(n)&=-\mu_{q,2}(n)-\mu_{q,3}(n)-\mu_{q,4}(n).
		\end{aligned}
	\end{equation}
	{\bf Acknowledgment}: The second author acknowledges the support from Project No. CRG/2019/00200/MS of the Science and Engineering Research Board, Department of Science and Technology, New Delhi, India.
	%%%%%%%%%%%%%%%%%%%%%%%%%%%%%%%%%%%%%%%%%%% BIBLIOGRAPHY %%%%%%%%%%%%%%%%%%%%%%%%%%%%%%%%%%%%%%%%%%%%%%%%

\end{document}